# Tangent and Supporting Lines, Envelopes, and Dual Curves

Steven J. Kilner
Department of Mathematics
Monroe Community College
Rochester, New York 14623 USA

David L. Farnsworth
School of Mathematical Sciences
Rochester Institute of Technology
Rochester, New York 14623 USA

**Abstract**

A differentiable curve $y = y(x)$ is determined by its tangent lines and is said to be the envelope of its tangent lines. The coefficients of the curve's tangent lines form a curve in another space, called the dual space. There is a transformation between the original $x,y$-space and the dual space, such that points on the original curve and the curve of the coefficients of the tangent lines are transformed into each other. The dual space and the transformation depend upon the form that is used for the tangent lines. One choice is $y = mx + b$, so that the coordinates in the dual space are $m$ and $b$, where the curve representing the tangent lines is $b = b(m)$. Each point of the curve $b = b(m)$ in the dual $m,b$-space corresponds to a tangent line to the curve $y = y(x)$ in $x,y$-space. We present other choices of the form for the tangent lines and explore techniques for finding a curve from the tangent lines, transformations between the original space and a dual space and among dual spaces, and the differential equation, whose solutions are exactly the equations of the tangent lines and the original curve. Geometric constructions and other tools are used to find dual curves.

## CONTENTS

## 1. The Envelope from the Tangent Lines

Consider the problem of determining a curve from its tangent lines. We begin with an example, where the tangent lines to the curve are presented with three different parameterizations. The initial parameter has no apparent geometrical significance. Then, the parameter is taken to be the $x$-coordinate of the point of tangency. Finally, the parameter is the slope of the tangent lines. An equation for the curve is found from the tangent lines.

**Example 1.1** (*three parametric representations of the tangent lines to a curve*)**:** Figure 1.1 shows some tangent lines to a curve. The lines are parameterized by $k > 1$ and have the equation
$$y = (k-1)^3 x + 1/(k-1)^3. \tag{1.1}$$
The goal is to find the function $y = y(x)$ that has these tangent lines.

One method, which is justified in Theorem 1.2, is to apply $d/dk$ to (1.1) and solve for $k = k(x)$, which is substituted into (1.1). Applying $d/dk$ to (1.1) yields
$$0 = 3(k-1)^2 x - 3(k-1)^{-4},$$
so that
$$x = (k-1)^{-6} \tag{1.2}$$
and
$$k = x^{-1/6} + 1. \tag{1.3}$$
Substituting (1.3) into (1.1) gives
$$y = ((x^{-1/6} + 1) - 1)^3 x + 1/((x^{-1/6} + 1) - 1)^3$$
or
$$y = 2x^{1/2}, \tag{1.4}$$
whose tangent lines are (1.1). Equation (1.2) and $k > 1$ imply that $x > 0$. Figure 1.2 displays a graph of (1.4).

Alternatively, the tangent lines to (1.4) can be parameterized by the $x$-coordinate $h$ of the point of tangency as
$$y = y'(h)(x - h) + y(h) = h^{-1/2}(x - h) + 2h^{1/2} = h^{-1/2}x - h^{1/2} + 2h^{1/2} = h^{-1/2}x + h^{1/2}. \tag{1.5}$$
To obtain (1.1), replace the parameter $h$ with $(k-1)^{-6}$, which gives
$$y = ((k-1)^{-6})^{-1/2}x + ((k-1)^{-6})^{1/2} = (k-1)^3 x + (k-1)^{-3}.$$

A third choice for the parameter of the tangent lines is their slope $m$. For the family in (1.1), the coefficient of $x$ is the slope $m$ in the standard form
$$y = mx + b. \tag{1.6}$$

Setting
$$m(k) = (k-1)^3 \qquad (1.7)$$
with $m > 0$ and solving (1.7) for $k$ gives
$$k = m^{1/3} + 1. \qquad (1.8)$$
Substituting (1.7) and (1.8) into (1.1) supplies
$$y = mx + 1/((m^{1/3} + 1) - 1)^3 = mx + 1/m, \qquad (1.9)$$
which uses $m$ to parameterize the set of tangent lines to (1.4), $y = 2x^{1/2}$ for $x > 0$.

Equation (1.9) shows the $y$-intercept $b = b(m) = 1/m$ as a function of slope $m$. The equations of the tangent lines of (1.4) are parameterized by $k$, $h$, and $m$ in (1.1), (1.5), and (1.9), respectively. Because $m$ has a geometrical significance and is a substitute for $dy/dx$, often it is used subsequently as the parameter.

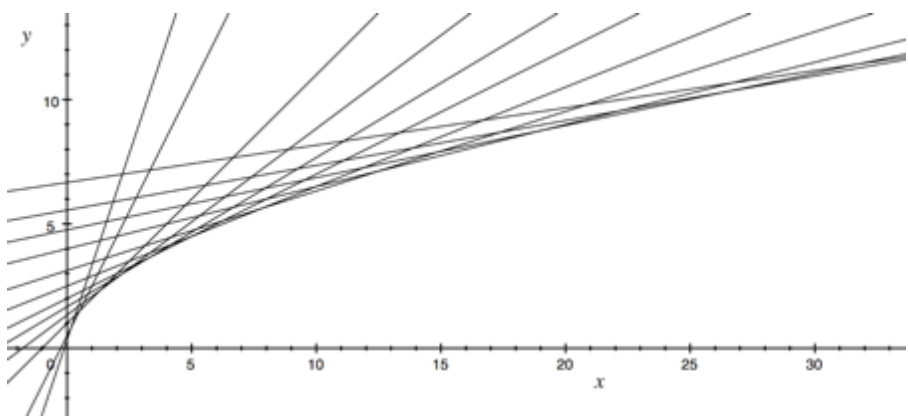

FIGURE 1.1: A set of tangent lines

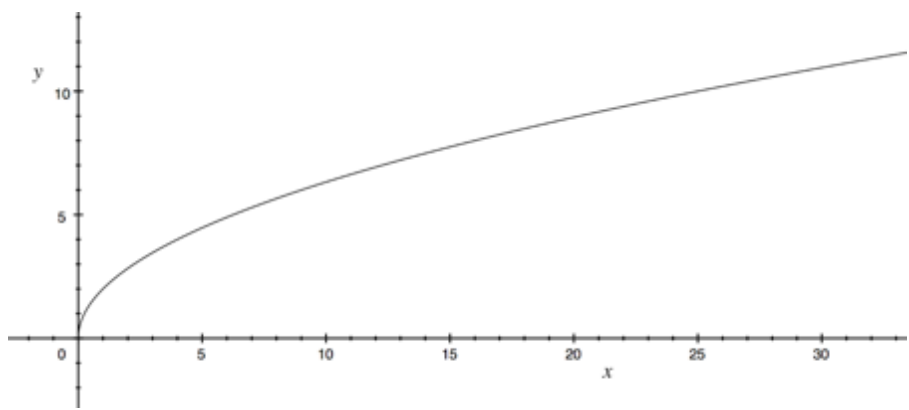

FIGURE 1.2: The curve $y = 2x^{1/2}$ for $x > 0$ whose tangent lines are in (1.1) and in Figure 1.1

In the plane, an *envelope* of a family of curves is a curve that is tangent to each member of the family. The envelope comprises those points of tangency. When the family of curves is composed of lines, the lines are the tangent lines to the envelope. The family of lines that is represented in (1.1), (1.5), and (1.9) have the curve (1.4) as their envelope.

The validity of the long-standing method of differentiating the equation of the tangent lines with respect to the parameter, solving for the parameter in terns of $x$, and substituting that function for the parameter into the equation of the tangent lines is established in Theorem 1.1 [28].

**Theorem 1.1** (*method of deriving the envelope from the tangent lines by differentiating with respect to the parameter*)**:** Consider the family of tangent lines

$$y = m(k)x + b(k), \qquad (1.10)$$

which is in the form based on the slope and the $y$-intercept. Applying the derivative $d/dk$ with respect to the parameter $k$ to (1.10) gives the $x$-coordinate as a function of $k$, and substituting that function into (1.10) gives the $y$-coordinate of a parametric representation of the envelope of the lines (1.10).

**Proof:** By changing $k$ by the small $\Delta k$, obtain the neighboring tangent line

$$y = m(k + \Delta k)x + b(k + \Delta k). \qquad (1.11)$$

The lines (1.10) and (1.11) intersect at a point near their points of tangency. As $\Delta k$ approaches zero, the point of intersection and the two points of tangency approach each other, and, in the limit, this results in a point on the curve.

The point of intersection is obtained by subtracting (1.10) from (1.11) giving

$$0 = (m(k + \Delta k) - m(k))x + b(k + \Delta k) - b(k). \qquad (1.12)$$

Dividing (1.12) by $\Delta k$ and letting $\Delta k \to 0$ gives

$$0 = m'(k)x + b'(k). \qquad (1.13)$$

The point of tangency that corresponds to the tangent line in (1.10) must satisfy equations (1.10) and (1.13). An explicit relationship between $x$ and $y$ on the curve can be obtained, if it is possible to eliminate the parameter $k$ between these two equations.

Assuming $m'(k) \neq 0$, (1.13) gives

$$x = -b'(k)/m'(k). \qquad (1.14)$$

Inserting (1.14) into (1.10) gives

$$y = m(k)(-b'(k)/m'(k)) + b(k). \qquad (1.15)$$

Equations (1.14) and (1.15) are a parametric representation of the envelope of the tangent lines (1.10). ∎

If (1.14) can be solved for $k$ as the function $k(x)$, then $k$ can be replaced in (1.15) to give the envelope as a function of $x$. For example, (1.1) and (1.3) yield (1.4) in Example 1.1.

Theorem 1.2 considers the case in Theorem 1.1 in which the parameter in the family of tangent lines is the $x$-coordinate of the point of contact of each tangent line with the envelope. Replacing that parameter with $x$ in the equation of the tangent lines provides the equation of the envelope. In Example 1.1, replacing $h$ by $x$ in (1.5) gives the envelope (1.4).

**Theorem 1.2** (*special case of Theorem 1.1 in which k is the x-coordinate*)**:** If the parameter $k$ in (1.10) is the $x$-coordinate of the point of tangency of the envelope, then the envelope is obtained by substituting $x$ for $k$ in (1.10).

**Proof:** If the parameter $k$ is the $x$-coordinate, then (1.10) can be written
$$y = m(k)(x - k) + y(k).$$
Applying $d/dk$ and recalling that $y'(k) = m(k)$ gives
$$0 = m'(k)x - m'(k)k - m(k) + y'(k) \text{ and } k = x. \quad \blacksquare$$

## 2. A Dual Space

Lines in $x,y$-space are determined by their slope $m$ and $y$-intercept $b$ in the slope-intercept form (1.6), with the exception of vertical lines, which do not have a finite slope. The two-dimensional space with coordinates $m$ and $b$ is called a *dual space* to $x,y$-space [8, p. 308–310; 31]. Each point in $m,b$-space identifies a line with equation $y = mx + b$ [29].

Consider the smooth curve $y = y(x)$ and its tangent lines. Each non-vertical tangent line is identified by a pair $(m,b)$ in the dual space. In Example 1.1, the tangent lines of $y = 2x^{1/2}$ in (1.4) are given by lines (1.9). The curve
$$b = b(m) = 1/m \tag{2.1}$$
for $m > 0$ in $b,m$-space contains points that provide the coefficients of the tangent lines. Figure 2.1 is a graph of (2.1).

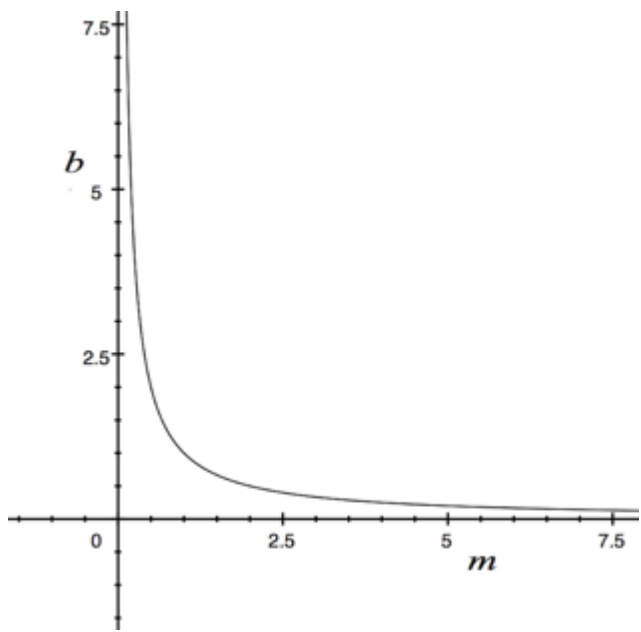

FIGURE 2.1: The dual curve $b = b(m) = 1/m$ with $m > 0$ in the dual $m,b$-space for $y = y(x) = 2x^{1/2}$ with $x > 0$ in $x,y$-space

A different dual space is introduced in Section 3 and used in Sections 4–6, 8.2, and 9. One reason for this replacement is that, if $y = y(x)$ is concave down (up), then the dual curve in $m,b$-space is concave up (down), as shown in Figures 1.2 and 2.1. For twice differentiable functions, it can be shown by using the method of proof of Theorem 3.1 that the signs of the second derivatives $y''(x)$ and $b''(m)$ have opposite signs at corresponding points.

Another, much more important reason is that a desirable property of the transformation between the original space and the dual space is that it be reflexive or involutive, that is, the dual of the dual of a curve would be the original curve. However, the dual in $m,b$-space of $y = y(x) = 1/x$ for $x > 0$ from (2.1) is $b = b(m) = 2(-m)^{1/2}$ for $m < 0$, which is not (4), $y = 2x^{1/2}$ for $x > 0$, even with the variables renamed.

A lack of symmetry in the definition for the dual $m,b$-space highlights the problem of the minus signs in the concavity and missing reflexivity by a negation in the example. For $m,b$-space, the tangent lines are written as (1.6), which can be re-expressed $y - b = mx$. The minus sign on the left-hand side displays the breaking of symmetry between the pairs $(x,y)$ and $(m,b)$.

## 3. An Alternative Dual Space and the Legendre Transformation

The property of reflexivity and the preservation of concavity can be achieved by using the slope $m$ and the negation of the $y$-intercept. This standard form for the tangent lines is

$$y = mx - d. \quad (3.1)$$

Theorems 1.1 and 1.2 apply by replacing $b$ by $-d$.

In Example 1.1,

$$y = y(x) = 2x^{1/2} \quad (3.2)$$

for $x > 0$ has tangent lines

$$y = mx - (-1/m) \text{ for } m > 0$$

from (1.9). The dual curve to (3.2) in $m,d$-space is

$$d = d(m) = -1/m. \quad (3.3)$$

The original curve is concave down, because $y''(x) = (-1/2)x^{-3/2} < 0$. The dual curve is also concave down, because $d''(m) = -2m^{-3} < 0$. The maintenance of the sign of the second derivative is shown in Theorem 3.1.

Further for Example 1.1, the tangent line to

$$y = y(x) = -1/x$$

from (3.3) at $x = h$ is

$$y = h^{-2}x - 2h^{-1}. \quad (3.4)$$

Because $m = h^{-2}$, $h = m^{-1/2}$ and (3.4) becomes

$$y = mx - (2m^{1/2}),$$

so that the dual curve in $m,d$-space is

$$d = d(m) = 2m^{1/2}$$

for $m > 0$, which is (3.2) with the necessary renaming of the variables. The reflexivity or involutive property for $x,y$-space and $m,d$-space is shown in Theorem 3.2.

The function $d(m)$ is called the *Legendre transformation of $y(x)$* and is written

$$\mathcal{L}\{y(x)\}(m) = d(m).$$

The dual $m,d$-space and the Legendre transformation are the focus of subsequent sections.

Maintenance of convexity and concavity by the Legendre transformation is shown in Theorem 3.1 [24; 34, pp. 257–261; 43], and reflexivity is shown in Theorem 3.2 [18, pp. 24–25; 27, p. 143; 34, pp. 104, 258].

**Theorem 3.1** (*preservation of strict convexity and concavity*)**:** If $\mathcal{L}\{y(x)\}(m) = d(m)$, then, where both functions possess second derivatives,
$$d''(m) = 1/y''(x). \tag{3.5}$$
**Proof:** Applying $d/dx$ to (3.1) and recalling that $m = y'(x)$ gives
$$y'(x) = y'(x) + xy''(x) - d'(m)y''(x).$$
For $y''(x) \neq 0$,
$$d'(m) = x. \tag{3.6}$$
Applying $d/dx$ gives
$$d''(m)y''(x) = 1,$$
which is (3.5). ∎

Because of the reflexive property, the slopes and negation of the $d$-intercepts of lines in $m,d$-space are given by the coordinates of the points of the corresponding curve $y = y(x)$ in $x,y$-space in the same manner as the slopes and negation of the $y$-intercepts of lines in $x,y$-space are given by the coordinates of the points of the corresponding curve $d = d(m)$ in $m,d$-space. Reflexivity assures that properties of the Legendre transformation are shared by the inverse Legendre transformation.

**Theorem 3.2** (*reflective or involutive property*)**:** After making any necessary changes in the names of the variables,
$$\mathcal{L}\{\mathcal{L}\{y(x)\}(m)\}(x) = y(x).$$
**Proof:** From (3.1), $y(x) = mx - d(m)$, where the variable is $x$. Solving (3.1) for $d$ and substituting $x = y'^{-1}(m)$ from $m = y'(x)$ gives
$$\mathcal{L}\{y(x)\}(m) = d(m) = mx - y(x) = my'^{-1}(m) - y(y'^{-1}(m)), \tag{3.7}$$
which is a formula for the Legendre transformation in any region where the function possesses an invertible derivative. Using (3.7) for the Legendre transformation of $d(m)$ into $x,y$-space gives
$$\mathcal{L}\{\mathcal{L}\{y(x)\}(m)\}(x) = \mathcal{L}\{d(m)\}(x) = xd'^{-1}(x) - d(d'^{-1}(x)).$$
From (3.6), $x = d'(m)$, so that
$$\mathcal{L}\{\mathcal{L}\{y(x)\}(m)\}(x) = \mathcal{L}\{d(m)\}(x) = xd'^{-1}(d'(m)) - d(d'^{-1}(d'(m))) = xm - d = y(x).$$
If the derivative is invertible on subdomains, the method of proof may be applied to each subdomain. ∎

## 4. Three Examples

Three examples are explored below for the dual space $m,d$-space. Examples 4.2 and 4.3 illustrate how dual curves, not just dual functions, can be analyzed.

**Example 4.1** (*the dual curves $y(x) = x^p/p$ and $d(m) = m^q/q$ for $1/p + 1/q = 1$*): Consider the curve
$$y = x^p/p, \tag{4.1}$$
with any real value of $p$ except 0 and 1 and $x > 0$, if necessary in order to maintain all real variables. We show that the dual curve is
$$d = m^q/q, \tag{4.2}$$
where $q$ determined by
$$1/p + 1/q = 1. \tag{4.3}$$

Using the form (3.1), that is, $y = mx - d$, for the tangent lines, the tangent line to (4.1) at $x = h$ is
$$y = h^{p-1}(x - h) - (-h^p/p) = h^{p-1}x - ((p-1)/p)h^p.$$
The slope is
$$m(h) = h^{p-1}, \tag{4.4}$$
and the negation of the $y$-intercept is
$$d = ((p-1)/p)h^p. \tag{4.5}$$
From (4.4),
$$h = m^{1/(p-1)}. \tag{4.6}$$
Substituting (4.6) into (4.5) gives the dual curve
$$d(m) = ((p-1)/p)m^{p/(p-1)}. \tag{4.7}$$
From (4.3),
$$1/q = 1 - 1/p = (p-1)/p \tag{4.8}$$
and
$$q = p/(p-1). \tag{4.9}$$
Substituting (4.8) and (4.9) into (4.7) gives the dual curve (4.2) in $m,d$-space.

Example 1.1 is the case of (4.1) and (4.2) where $p = \frac{1}{2}$ and $q = -1$.

**Example 4.2** (*the dual curves $y(x) = x^3 - x$ and $d(m) = 2((m + 1)/3)^{3/2}$*): Figure 4.1 contains tangent lines to $y = x^3 - x$, and Figure 4.2 contains the envelope of the tangent lines, that is, $y = x^3 - x$. The tangent lines are
$$y = (3k^2 - 1)(x - k) + (k^3 - k) = (3k^2 - 1)x - 2k^3.$$
Setting $m = (3k^2 - 1)$ gives $k^2 = (m + 1)/3$, and the tangent lines become
$$y = mx - 2((m + 1)/3)^{3/2}.$$
The dual curve is
$$d(m) = 2((m + 1)/3)^{3/2}, \tag{4.10}$$
where both the positive and the negative square roots are taken.

The dual curve appears in Figure 4.3. The minimum of $m$ is $m(0) = -1$, and $m = dy/dx \geq -1$. Corresponding points $A$, $B$, and $C$ are labeled in Figures 4.2 and 4.3 in order to indicate the directions "followed" by the two curves' corresponding points. The point of inflection at $(x,y) = (0,0)$ and the cusp at $(m,d) = (-1, 0)$ in Figure 4.2 and 4.3 are the locations of changes in concavity.

As Figure 4.1 shows, some pairs of points of the curve in $x,y$-space share slopes, but the negations of their $y$-intercepts are different. The dual curve $d = d(m)$ is not a function, but it describes a curve in $m,d$-space.

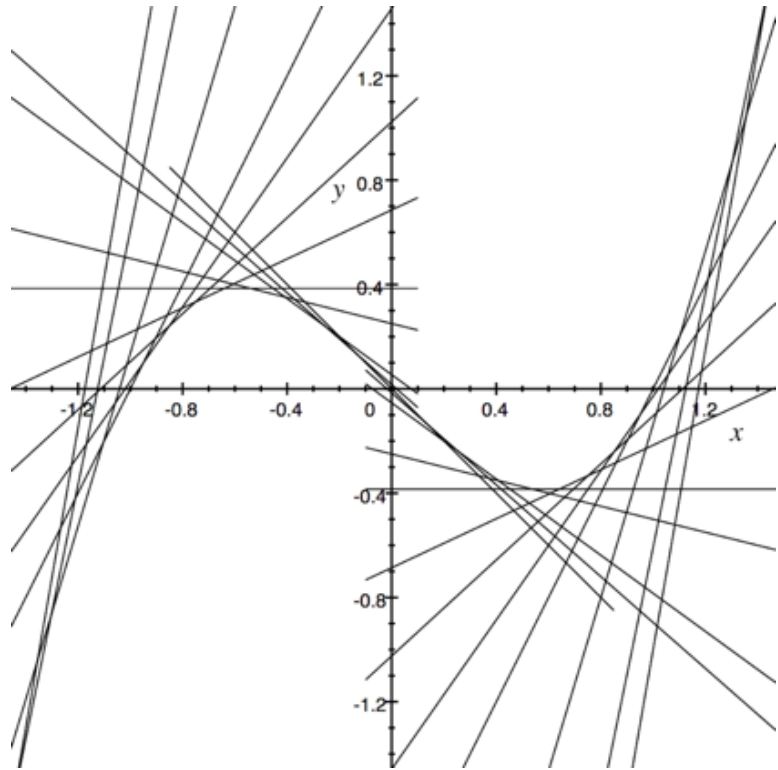

FIGURE 4.1: Tangent lines to $y = x^3 - x$

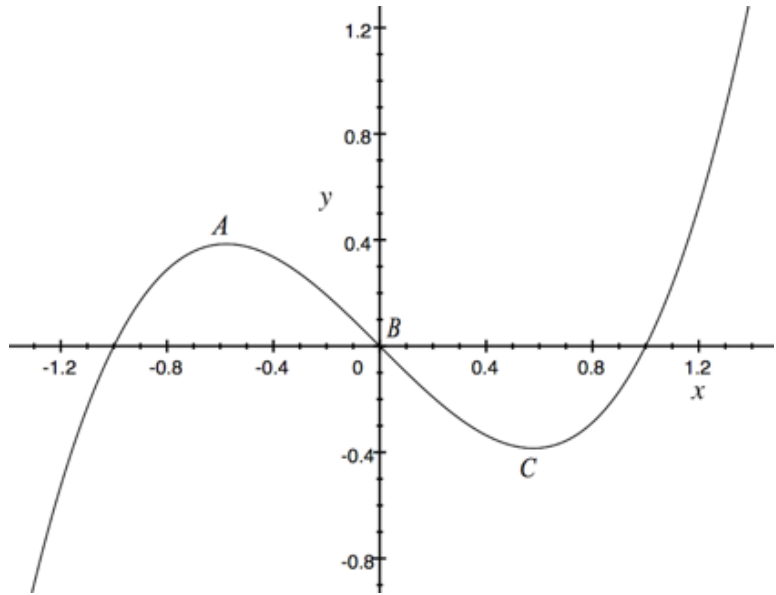

FIGURE 4.2: Envelope $y = x^3 - x$ of the tangent lines in Figure 4.1

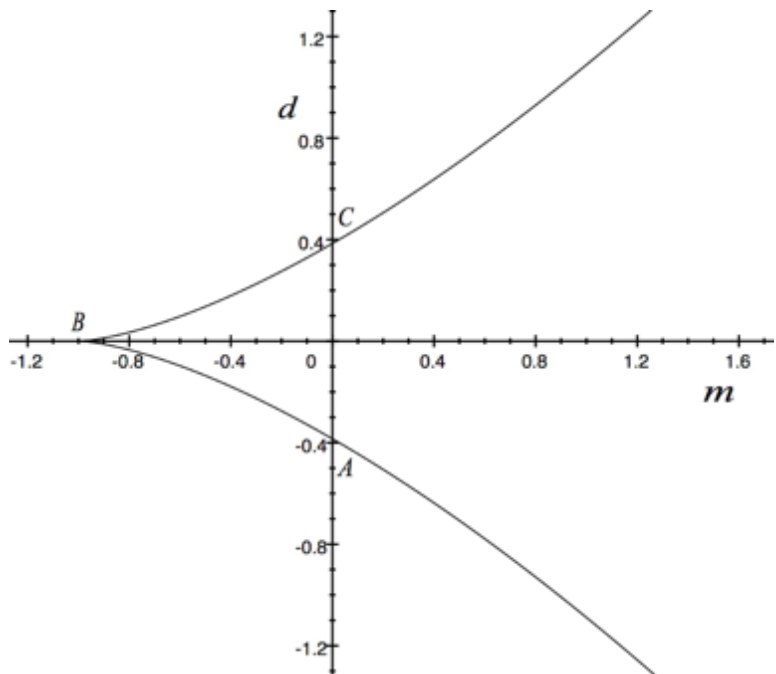

FIGURE 4.3: Dual curve $d(m) = 2((m + 1)/3)^{3/2}$ to $y = x^3 - x$

**Example 4.3** (*a curve with a shared tangent line at two points and three vertical tangent lines*): Consider the curve

$$y(x) = \begin{cases} \sqrt{1-(x+1)^2} & \text{for } -2 \le x \le 0 \\ \sqrt{1-(x-1)^2} & \text{for } 0 < x \le 2 \end{cases}, \qquad (4.11)$$

which is composed of two semicircles and appears in Figure 4.4. The points $(-1,1)$ and $(1,1)$ share the tangent line $y = 1$, and, hence, have the same slope $m = 0$ and negations of the $y$-intercept $d = -1$. The dual points of their tangent lines coincide at $(m,d) = (0,-1)$. Three points have vertical tangent lines. As $x$ approaches $\pm 2$ at the left and right boundary points, the slope and $y$-intercepts become nonfinite, so the tangent lines would correspond to points at infinity in the dual space. The internal tangent line at the origin is at a cusp of $y = y(x)$. It would correspond to a point at infinity with $d = 0$.

The dual curve $d = d(m) = \mathcal{L}\{y(x)\}(m)$ to the curve in Figure 4.4 is

$$d(m) = \begin{cases} \sqrt{1+m^2} - m \\ \sqrt{1+m^2} + m \end{cases}, \qquad (4.12)$$

which is defined for all values of $m$ and appears in Figure 4.5.

In order to indicate the paths of corresponding points on the two curves in Figures 4.4 and 4.5, corresponding points are labeled. As $x$ increases, in Figure 4.4, the two semicircular arcs are traced and the points $A$ through $G$ are met in alphabetical order. In Figure 4.5, the self-intersection or crossing point $(m,d) = (0,-1)$ is dual of the tangent line $y = 1$, which is shared by the two points $(x,y) = (-1,1)$ and $(1,1)$ in Figure 4.4.

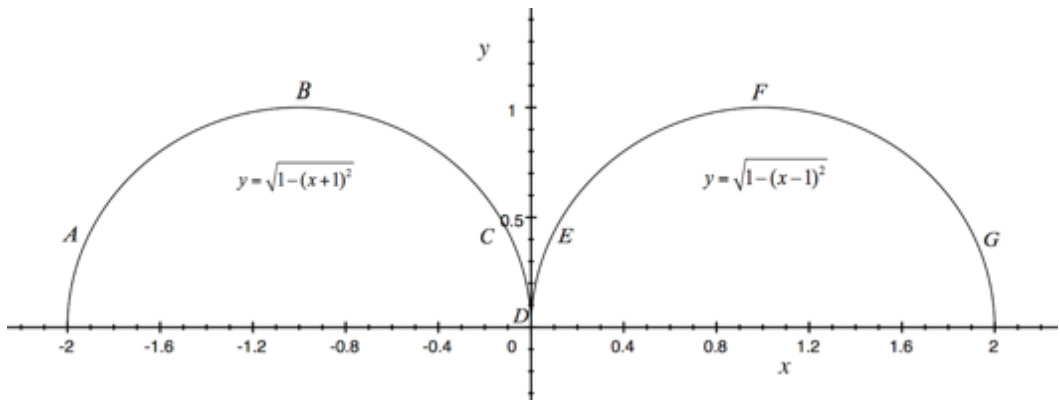

FIGURE 4.4: A continuous curve in (4.11) with the shared tangent line $y = 1$ by two points and three vertical tangent lines

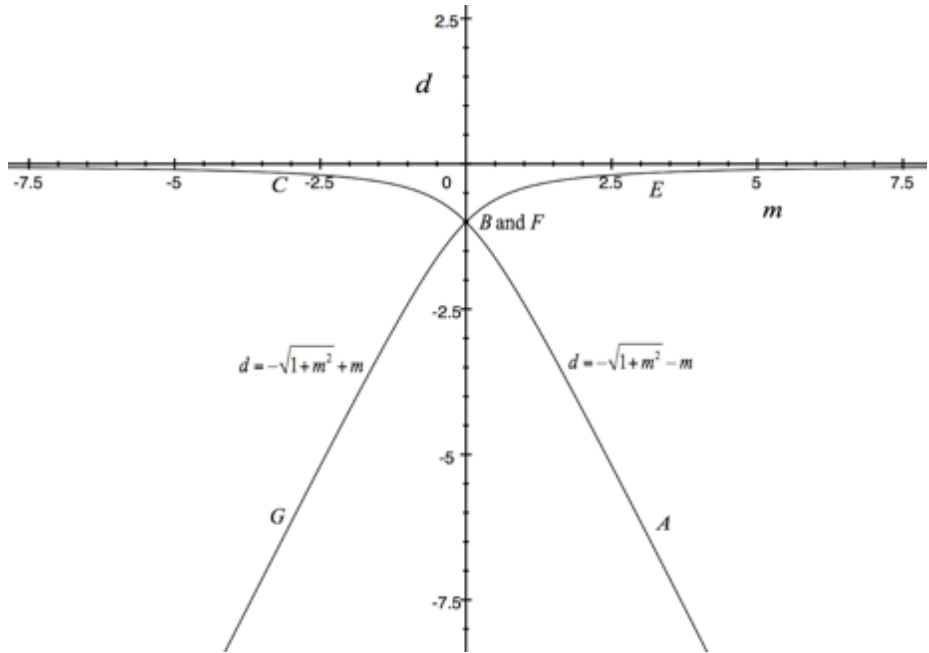

FIGURE 4.5: Curve in (4.12), which is dual to the curve in Figure 4.4

## 5. Properties of the Legendre Transformation

Theorem 5.1 supplies another way to find a dual curve up to an additive constant [43].

**Theorem 5.1** ($d'(y'(x)) = x$ and $y'(d'(m)) = m$): The first derivatives of the pair of differentiable functions

$$d(m) = \mathcal{L}\{y(x)\}(m) \text{ and } y(x) = \mathcal{L}\{d(m)\}(x)$$

are inverses, that is,

$$d'(y'(x)) = x \tag{5.1}$$

and

$$y'(d'(m)) = m. \tag{5.2}$$

**Proof:** From (3.6),

$$x = d'(m). \tag{5.3}$$

Also,

$$m = y'(x). \tag{5.4}$$

Substituting $m$ from (5.4) into the right-hand side of (5.3) gives (5.1). Substituting $x$ from (5.3) into the right-hand side of (5.4) gives (5.2). ∎

Most transformations have a representation as an integral. The Legendre transformation has one, which is displayed in Corollary 5.1. Refer to Theorem 9.1, as well.

**Corollary 5.1** (*representation of the Legendre transformation as an integral*): The Legendre transformation of the function $y = y(x)$ is

$$\mathcal{L}\{y(x)\}(m) = d(m) = \int y'^{-1}(m) dm \tag{5.5}$$

in domains where the function $y(x)$ possesses an invertible derivative.

**Proof:** Formula (5.5) is a direct consequence of (5.2), which implies that

$$d'(m) = y'^{-1}(m)$$

and

$$\mathcal{L}\{y(x)\}(m) = d(m) = \int y'^{-1}(m) dm. \blacksquare$$

For an example of Corollary 5.1, consider $y = y(x) = x^p/p$, with $p > 1$ from Example 4.1. Then, $m = y'(x) = x^{p-1}$ and $y'^{-1}(m) = m^{1/(p-1)}$. From (5.5),

$$\mathcal{L}\{y(x)\}(m) = d(m) = \int y'^{-1}(m) dm = \int m^{1/(p-1)} dm = \frac{m^{1/(p-1)+1}}{1/(p-1)+1} = \frac{m^q}{q},$$

using the definition of $q$ in (4.3). The constant of integration is zero because $m(0) = 0$.

Table 5.1 contains pairs that are Legendre transformations of each other. Theorem 3.2 says that the table can be read right-to-left and left-to-right.

TABLE 5.1: Legendre-transformation pairs, where $\mathcal{L}\{y(x)\}(m) = d(m)$

| | **g(x)** | **$\mathcal{L}\{g(x)\}(m)$** |
|---|---|---|
| 1. | $ay(x)$ | $ad(m/a)$, $a \neq 0$ |
| 2. | $y(ax)$ | $d(m/a)$, $a \neq 0$ |
| 3. | $y(x) + a$ | $d(m) - a$ |
| 4. | $y(x + a)$ | $d(m) - am$ |
| 5. | $y^{-1}(x)$ | $-md(1/m)$, $m \neq 0$ |
| 6. | $cy(sx + t) + bx + a$ | $cd((m-b)/(cs)) - t(m-b)/s - a$, $c \neq 0$, $s \neq 0$ |
| 7. | $x^p/p$ | $m^q/q$, $x > 0$, $m > 0$, $p > 1$ $1/p + 1/q = 1$ |
| 8. | $\sqrt{1-x^2}$ | $-\sqrt{1+m^2}$, $-1 < x < 1$ |
| 9. | $(1-x^p)^{1/p}$ | $-(1+(-m)^q)^{1/q}$, $0 \leq x < 1$, $m < 0$, $p > 1$, $1/p + 1/q = 1$ |
| 10. | $\sqrt{x^2 - 1}$ | $\sqrt{m^2 - 1}$, $|x| \geq 1$, $|m| \geq 1$ |
| 11. | $e^x$ | $m \ln m - m$, $m > 0$ |
| 12. | $\ln x$ | $1 + \ln m$, $x > 0$, $m > 0$ |
| 13. | $x \ln x$ | $e^{m-1}$, $x > 0$ |
| 14. | $\cos x$ | $-m \sin^{-1} m - \sqrt{1-m^2}$, $0 \leq x \leq \pi/2$, $-1 \leq m \leq 0$ |

## 5.1 Entries 1–5 of Table 5.1

Entries 1–5 of Table 5.1 can be seen for regions over which the curves are sufficiently differentiable as follows. From the definition of $\mathcal{L}\{y(x)\}(m) = d(m)$, we have the formula

$$\mathcal{L}\{y(x)\}(m) = d(m) = m{y'}^{-1}(m) - y({y'}^{-1}(m)), \qquad (3.7)$$

which gives $d(m)$ in terms of only the function $y(x)$ along with its derivative and its inverse.

For Entry 1, consider $g(x) = ay(x)$, whose tangent line at $x = c$ is

$$y = ay'(c)x + (ay(c) - cay'(c)).$$

Setting $m = ay'(c)$, so that $c = {y'}^{-1}(m/a)$, gives

$$\begin{aligned}
\mathcal{L}\{ay(x)\}(m) &= cay'(c) - ay(c) \\
&= a{y'}^{-1}(m/a)y'({y'}^{-1}(m/a)) - ay({y'}^{-1}(m/a)) \\
&= a(m/a){y'}^{-1}(m/a) - ay({y'}^{-1}(m/a)) \\
&= a((m/a){y'}^{-1}(m/a) - y({y'}^{-1}(m/a))) \\
&= ad(m/a)
\end{aligned}$$

for $d(m)$ in (3.7).

For Entry 2, consider $g(x) = y(ax)$, whose tangent line at $x = c$ is

$$y = ay'(ac)x + (y(ac) - cay'(ac)).$$

Setting $m = ay'(ac)$, so that $c = (1/a){y'}^{-1}(m/a)$, gives

$$\begin{aligned}
\mathcal{L}\{y(ax)\}(m) &= cay'(ac) - y(ac) \\
&= a(1/a){y'}^{-1}(m/a)y'({y'}^{-1}(m/a)) - y((a)(1/a){y'}^{-1}(m/a)) \\
&= (m/a){y'}^{-1}(m/a) - y({y'}^{-1}(m/a)) \\
&= d(m/a)
\end{aligned}$$

for $d(m)$ in (3.7).

For Entry 3, consider $g(x) = y(x) + a$, whose tangent line at $x = c$ is

$$y = y'(c)x + (y(c) + a - cy'(c)).$$

Setting $m = y'(c)$, so that $c = {y'}^{-1}(m)$, gives

$$\begin{aligned}
\mathcal{L}\{y(x) + a\}(m) &= cy'(c) - y(c) - a \\
&= {y'}^{-1}(m)y'({y'}^{-1}(m)) - y({y'}^{-1}(m)) - a \\
&= m{y'}^{-1}(m) - y({y'}^{-1}(m)) - a \\
&= d(m) - a
\end{aligned}$$

for $d(m)$ in (3.7).

For Entry 4, consider $g(x) = y(x + a)$, whose tangent line at $x = c$ is
$$y = y'(c + a)x + (y(c + a) - cy'(c + a)).$$
Setting $m = y'(c + a)$, so that $c = y'^{-1}(m) - a$, gives
$$\mathcal{L}\{y(x + a)\}(m) = cy'(c + a) - y(c + a)$$
$$= (y'^{-1}(m) - a)y'(y'^{-1}(m) - a + a) - y(y'^{-1}(m) - a + a)$$
$$= m(y'^{-1}(m) - a) - y(y'^{-1}(m))$$
$$= d(m) - am$$
for $d(m)$ in (3.7).

For Entry 5, consider $g(x) = y^{-1}(x)$, whose tangent line at $x = c$ is
$$y = y^{-1\prime}(c)x + (y^{-1}(c) - cy^{-1\prime}(c)).$$
The Legendre transformation is
$$\mathcal{L}\{y^{-1}(x)\}(m) = cy^{-1\prime}(c) - y^{-1}(c).$$
Designate $a = y^{-1}(c)$. Then, $c = y(a)$. Set $m = y^{-1\prime}(c)$, which is equal to $1/y'(a)$ [10, pp. 219–221; 12, pp. 144–147; 25, pp. 345–346], so that $y'(a) = 1/m$. Then, $a = y'^{-1}(1/m) = y^{-1}(c)$ and $c = y(a) = y(y'^{-1}(1/m))$. The Legendre transformation is
$$\mathcal{L}\{y^{-1}(x)\}(m) = y(y'^{-1}(1/m)) m - y'^{-1}(1/m)$$
$$= -m((1/m)y'^{-1}(1/m) - y(y'^{-1}(1/m))) = -md(1/m)$$
for $d(m)$ in (3.7).

## 5.2 Other Entries of Table 5.1

Entry 6 is found by using Entries 1–4 reading both right-to-left and left-to-right. Entry 7 is Example 4.1. In Example 4.2, the transformation $\mathcal{L}\{x^3 - x\}(m) = 2((m + 1)/3)^{3/2}$ can be found by using Entry 6 for $a = 0$, $b = -1$, $c = 3$, $s = 1$, and $t = 0$ and Entry 7 with $p = 3$, $q = 3/2$. The inverse Legendre transformation $\mathcal{L}\{2((m + 1)/3)^{3/2}\}(x) = x^3 - x$ is obtained by reading Entries 6 and 7 right-to-left.

For Entry 11, the tangent line to $y = e^x$ at $x = c$ is
$$y = e^c x + e^c - ce^c.$$
Setting $m = e^c$ so that $c = \ln m$ gives
$$\mathcal{L}\{e^x\}(m) = ce^c - e^c = (\ln m)e^{\ln m} - e^{\ln m} = m\ln m - m.$$
Alternatively, for $y = e^x$, we have $y'(x) = e^x = m$. Substituting $y = e^x$ and $y'^{-1}(m) = \ln m$, into (3.7) gives
$$\mathcal{L}\{e^x\}(m) = m\ln m - e^{\ln m} = m\ln m - m.$$

Entry 12 can be found using Entries 5 and 11 as follows. For $y = e^x$, $x = \ln y$. Writing $y^{-1}(x) = \ln x$,
$$\mathcal{L}\{\ln x\}(m) = \mathcal{L}\{y^{-1}(x)\}(m) = -md(1/m) = -m((1/m)\ln(1/m) - 1/m) = 1 + \ln m.$$

For Entry 13, taking $a = -1$ in Entry 4 yields $\mathcal{L}\{y - 1\}(m) = d(m) + m$. Applying that to Entry 11 gives

$$\mathcal{L}\{e^{x-1}\}(m) = (m\ln m - m) + m = m\ln m$$

and by Theorem 3.2

$$\mathcal{L}\{x\ln x\}(m) = e^{m-1}.$$

Some functions have the same functional form as their dual functions. One example is Entry 7 with $p = q = 2$, and another is Entry 10. A third example is obtained by combining Entries 3 and 12 to yield

$$\mathcal{L}\{a + \ln x\}(m) = 1 - a + \ln m.$$

and selecting $a = 1/2$. The whole curves have the same functional form, but, except for Entry 7 with $p = q = 2$, the coordinates in $x,y$-space and the corresponding point in $m,d$-space do not have the same numerical values, with the possible exception of isolated points. In order for the curves to be in one-to-one correspondence with the same coordinates point by point, the slope must be equal to the $x$-coordinate in $x,y$-space. The only function satisfying that condition is $y = x^2/2 + C$, because its slope is $y'(x) = x$, but, according to Entry 3 of Table 5.1, $\mathcal{L}\{x^2/2 + C\}(m) = m^2/2 - C$ $\neq m^2/2$, unless $C = 0$. This is confirmed by Theorem 8.3 and is shown in Section 8.1 for this dual space.

## 6. Clairaut's Differential Equation and the Legendre Transformation

A family of tangent lines and its envelope are precisely the solutions to a differential equation, which is called Clairaut's differential equation. A function that appears in the differential equation is the Legendre transformation of the equation's singular solution.

### 6.1 Clairaut's Differential Equation

Theorems 6.1 and 6.2 state that both the family of tangent lines in the standard form (3.1), that is,

$$y = mx - d$$

and its envelope are the solutions of *Clairaut's differential equation*

$$y = xy' + f(y'). \tag{6.1}$$

The set of tangent lines is called the *general solution* of (6.1) and contains a parameter, which is a constant of integration. The envelope is the *singular solution* of (6.1) and is not included in the general solution for any value of the parameter [9, pp. 39–40; 33, pp. 281–284; 44, pp. 75–77]. In Section 6.2, we see that Equation (6.1) contains the Legendre transformation of the singular solution, because the function $f(m)$ is $-d(m)$, where $d = d(m)$ is the envelope's dual curve in $m,d$-space.

**Theorem 6.1** (*solutions of Clairaut's differential equation*): Consider the family of tangent lines
$$y = m(k)x - d(m(k))$$
to the curve $y = y(x)$, where the parameter $k$ need not have any particular geometric significance. The tangent lines and their envelope are solutions to Clairaut's differential equation (6.1), where $f$ is the function such that
$$d(m(k)) = -f(m(k)). \tag{6.2}$$
**Proof:** The tangent line at $(h, y(h))$ in the form (3.1) is
$$y = y'(h)(x - h) + y(h) = y'(h)x + y(h) - hy'(h) = m(k)x - d(m(k)). \tag{6.3}$$
Because (6.3) is true for an interval of $x$-values,
$$y'(h) = m(k) \tag{6.4}$$
and
$$y(h) - hy'(h) = -d(m(k)). \tag{6.5}$$
If there is a function $f$ such that $d$ can be represented by (6.2), then from (6.5)
$$y(h) - hy'(h) = f(m(k)). \tag{6.6}$$
Using (6.4) to substitute for $m(k)$, (6.6) becomes
$$y(h) - hy'(h) = f(y'(h)). \tag{6.7}$$
Replacing the arbitrary $h$ by $x$ in (6.7) gives Clairaut's differential equation (6.1). ∎

**Theorem 6.2** (*solving Clairaut's differentia equation*): Clairaut's differential equation (6.1) has the general solution
$$y = Cx - f(C), \tag{6.8}$$
where $C$ a constant of integration, which is the family of tangent lines to the singular solution. The singular solution is the envelope of the general solution and can be represented parametrically by
$$x = -f'(t) \text{ and } y = f(t) - tf'(t). \tag{6.9}$$
If $f''(t) \neq 0$, that is, if $f$ is not a linear function, then, the singular solution is not included in the general solution.

**Proof:** Applying $d/dx$ to (6.1) gives $y' = y''x + y' + f'(y')y''$, so that $y'' = 0$ or $x = -f'(y')$. The general solution arises from $y'' = 0$, which gives $y' = C$. Substituting that into (6.1) yields (6.8). Substituting $x = -f'(y')$ into (6.1) and naming the parameter $t$ yields (6.9). To show that the singular solution is not part of the general solution, find
$$\frac{dy}{dx} = \frac{\frac{dy}{dt}}{\frac{dx}{dt}} = \frac{f'(t) - f'(t) - tf''(t)}{-f''(t)} = t,$$
which is not a constant, that is, the general solution cannot be a line. ∎

The singular solution (6.9) of Clairaut's differential equation is obtained by setting $k = m = t$ and $b(t) = f(t)$ in the envelope (1.14) and (1.15). Also, (6.9) is obtained by setting $m = t$ and $d(t) = -f(t)$ in (3.1) and (3.6).

The relationship (3.1) between dual curves and Clairaut's differential equation (6.1) are examples of the integration by parts formula $\int u\,dv = uv - \int v\,du$. For (3.1), set $u = m = y'(x)$ and $v = x = d'(m)$.

## 6.2 Relationship between Clairaut's Equation and the Legendre Transformation

Because the dual curve is represented in Clairaut's differential equation with $f(y')$ and the singular solution is the dual of that curve, solving Clairaut's differential equation provides another method to find dual curves and entries for Table 5.1. Conversely, Entries 7–14 of Table 5.1 might be considered an integral table for Clairaut's differential equation [36].

**Example 6.1** (*solutions of Clairaut's differential equation supply Legendre transformations and vice versa*)**:** Consider Clairaut's differential equation

$$y = xy' - \exp y'. \tag{6.10}$$

From (6.8), the general solution is

$$y = Cx - \exp C. \tag{6.11}$$

From (6.9), the singular solution is

$$x = \exp t \tag{6.12}$$

with

$$y = -\exp t + t \exp t. \tag{6.13}$$

From (6.12), $t = \ln x$, so that (6.13) becomes

$$y = -x + x \ln x. \tag{6.14}$$

In summary, the lines (6.11) are the tangent lines to the envelope (6.14). The Legendre transformation of (6.14) is $d(m) = \exp m$ from the second term on the right-hand side of (6.10). Solving Clairaut's differential equation has established the validity of Entry 11 in Table 5.1, read right-to-left, which is permissible by Theorem 3.2.

Given $y = y(x)$, substitution into Clairaut's differential equation in the format

$$y = xy' - d(y')$$

yields the function $d(m) = \mathcal{L}\{y(x)\}(m)$ as follows. Setting $y = y(x) = e^x$ in the differential equation gives $e^x = xe^x - d(e^x)$. Solving for $d$, this becomes

$$d(e^x) = (x - 1)e^x.$$

Substituting $m = y'(x) = e^x$, so that $x = \ln m$, gives

$$d(m) = (\ln m - 1)m = \mathcal{L}\{e^x\}(m).$$

Entry 11 can be employed to solve the differential equation. Consider Clairaut's differential equation (6.10),

$$y = xy' - \exp y'$$

and Entry 11 of Table 5.1,

$$\mathcal{L}\{e^x\}(m) = m\ln m - m.$$

The singular solution to (6.10) is $y = x\ln x - x$ as in (6.14).

Using this relationship with Clairaut's differential equation, the differential equation can be used to derive Entries 1–5 of Table 5.1. For Entry 1, writing $g(x) = ay(x)$ and replacing $y(x)$ by $g(x)/a$ in Clairaut's differential equation (6.1) give

$$g(x)/a = xg'(x)/a + f(g'(x)/a).$$

Multiplying by $a$ yields

$$g(x) = xg'(x) + af(g'(x)/a). \tag{6.15}$$

Comparing (6.15) with (6.1) and using (6.2) give

$$\mathcal{L}\{ay(x)\}(m) = \mathcal{L}\{g(x)\}(m) = -af(m/a) = ad(m/a).$$

For Entry 2, evaluating Clairaut's differential equation (6.1) at $ax$ gives

$$y(ax) = (ax)y'(ax) + f(y'(ax)). \tag{6.16}$$

Recall that prime ( ′ ) indicates derivative with respect to the function's argument. Then, (6.16) becomes

$$y(ax) = (ax)(dy(ax)/dx)(1/a) + f((dy(ax)/dx)(1/a)). \tag{6.17}$$

Cancelling $a$ on the right-hand side of (6.17), comparing it to (6.1), and using (6.2) give

$$\mathcal{L}\{y(ax)\}(m) = -f(m/a) = d(m/a).$$

For Entry 3, writing $g(x) = y(x) + a$ and replacing $y(x)$ by $g(x) - a$ in (6.1) give

$$g(x) - a = xg'(x) + f(g'(x)),$$

so that

$$g(x) = xg'(x) + f(g'(x)) + a. \tag{6.18}$$

Comparing (6.18) with (6.1) and using (6.2) give

$$\mathcal{L}\{y(a) + a\}(m) = \mathcal{L}(g(x)) = -(f(m) + a) = d(m) - a.$$

For Entry 4, evaluating (6.1) at $x + a$ gives

$$y(x + a) = (x + a)y'(x + a) + f(y'(x + a))$$

and

$$y(x + a) = (x + a)dy(x + a)/dx + f(dy(x + a)/dx),$$

which becomes

$$y(x + a) = xdy(x + a)/dx + f(dy(x + a)/dx) + ady(x + a)/dx \tag{6.19}$$

Comparing (6.19) to (6.1) and using (6.2) give

$$\mathcal{L}\{y(x + a)\}(m) = -(f(m) + am) = d(m) - am.$$

For Entry 5, evaluating (6.1) at $y^{-1}(x)$ gives

$$y(y^{-1}(x)) = y^{-1}(x)\, y'(y^{-1}(x)) + f(y'(y^{-1}(x))). \tag{6.19}$$

Because

$$y'(y^{-1}(x)) = 1/(y^{-1})'(x)$$

[10, pp. 219–221; 12, pp. 144–147; 25, pp. 345–346], equation (6.20) is

$$x = y^{-1}(x)\,(1/(y^{-1})'(x)) + f(1/(y^{-1})'(x)). \tag{6.21}$$

Multiplying (6.21) by $(y^{-1})'(x)$ and rearranging terms result in
$$y^{-1}(x) = x(y^{-1})'(x) + (y^{-1})'(x) f(1/1/(y^{-1})'(x)). \qquad (6.22)$$
Comparing (6.22) to (6.1) and using (6.2) give
$$\mathcal{L}\{y^{-1}(x)\}(m) = -m(f(1/m)) = md(1/m).$$

## 7. Choices for Dual Spaces and Relationships among Them

For each curve $y = y(x)$, which is the envelope in $x,y$-space, there is a corresponding curve in the dual space that expresses the relationship between the two coefficients in the equations of the family of the supporting or tangent lines. A feature of all dual spaces is that they preserve information [14]. Dual curves and the original curves can be created from each other. Even in Example 4.3, where the two points $(-1,1)$ and $(1,1)$ in $x,y$-space share a tangent line, there is no ambiguity about the identity of the dual curves. Each selection of a standard form for the tangent lines to a curve gives a different dual space. Five standard canonical forms for equations of lines are in Table 7.1. The choice among them depends upon the coefficients of interest, relevant applications, and importance of properties, such as reflexivity or maintenance of convexity.

TABLE 7.1: Standard formulas of lines

| Name | Equation | Parameters | Breakdown |
|---|---|---|---|
| 1. Slope-intercept | $y = mx + b$ | Slope $m$ and y-intercept $b$ | Vertical lines |
| 2. Slope-negation-of-the-$y$-intercept | $y = mx - d$ | Slope $m$ and negation of the $y$-intercept ($d$) | Vertical lines |
| 3. Dot-product or inverses of intercepts | $ux + vy = 1$ | Coordinate intercepts $1/u$ and $1/v$ | Lines containing $(0,0)$ |
| 4. Fixed point [17] | $y = mx + (1-m)t$ | Slope $m$ and $t$ such that $y(t) = t$ | Lines parallel to $y = x$ |
| 5. Polar | $(\cos\phi)x + (\sin\phi)y = D$ | Perpendicular distance $D > 0$ from $(0,0)$ to the line and counterclockwise angle $\phi$ from the positive $x$-axis to the direction in which $D$ is measured with $0 \le \phi < 2\pi$; lines through $(0,0)$ have coordinates $(0,\phi)$ with $0 \le \phi < \pi$, which are obtained from $(D,\phi)$ with $D \to 0^+$ | |

For any one-parameter family of lines, such as all the tangent lines to a single curve, both coefficients can be expressed as functions of a single parameter or else one parameter or coefficient can be presented as a function of the other. For standard forms 1–4 in Table 7.1, the customary relationships are $b = b(m)$, $d = d(m)$, $v = v(u)$ or, less often, $u = u(v)$, and $t = t(m)$, respectively.

There are simple relationships among the parameters or coefficients of the forms for lines. For example, relationships between form 1 and each of the Forms 3–5 are as follows:

- The slope-intercept coefficients $m$ and $b$ can be obtained from the dot-product coefficients $u$ and $v$ by $m = -u/v$ and $b = 1/v$. The inverse relationship between these two pairs of coefficients is $u = -m/b$ and $v = 1/b$.

- The slope-intercept coefficients $m$ and $b$ can be obtained from the fixed-point coefficients $m$ and $t$ by $m = m$ and $b = (1 - m)t$. The inverse relationship between these two pairs of coefficients is $m = m$ and $t = b/(1 - m)$.

- The slope-intercept coefficients $m$ and $b$ can be obtained from the polar parameters $D$ and $\phi$ by $m = -\cot \phi$ and $b = D/\sin \phi$. The inverse relationship between these two pairs of coefficients is $D = |b|/\sqrt{1 + m^2}$ and $\phi = \cot^{-1}(-m)$.

Form 3, $ux + vy = 1$, is used in vectors spaces [18, p. 23; 26, pp. 73–90; 38, pp. 1–5]. It is called the dot-product form because $ux + vy = (u,v)\bullet(x,y)^t = 1$. Also, $(1/u, 0)$ and $(0, 1/v)$ are the points of intersection of the line with the axes. For this form, differentiation can be employed to find the envelope, as in Theorem 1.1, and $x,y$-space and $u,v$-space are reflexive, as in Theorem 3.2.

**Example 7.1** (*change of parameters between different dual spaces*)**:** Entries in Table 5.1 can aid in finding dual curves in other spaces. For example, from the first bulleted item above and because the negation of the $y$-intercept $(-b)$ is $d$, the slope-negation-of-the-$y$-intercept coefficients are $m = -u/v$ and $d = -1/v$. Substituting those into Entry 9 of Table 5.1

$$d = -(1 + (-m)^q)^{1/q}$$

gives

$$-1/v = -(1 + (-(-u/v))^q)^{1/q}$$

or

$$1/v^q = (u^q + v^q)/v^q.$$

The dual curve to $x^p + y^p = 1$ or $y = (1 - x^p)^{1/p}$ in $x,y$-space is $d(m) = -(1 + (-m)^q)^{1/q}$ in $m,d$-space and is

$$v = (1 - u^q)^{1/q}$$

or $u^q + v^q = 1$ in $u,v$-space.

**Example 7.2** ($y = x\ln x - x$ *and two of its dual curves*)**:** The function $y = x\ln x - x$, two of its dual curves, and transformations among them are displayed in Figures 7.1–7.4. See Entry 11 of Table 5.1, read right-to-left. At each vertex of the triangle in Figure 7.1, there is

- the name of the space with its variables in the order of independent variable then dependent variable,

- the equation of a curve that corresponds to the curves at the other corners by the transformations displayed along the arrows,

- a generic tangent line to the curve that is expressed first with the value of the independent variable at the point of tangency as the parameter, second with the slope at the point of tangency as the coefficient of the independent variable and the parameter, and third in the dot-product form, and
- Clairaut's differential equation that is satisfied by both the curve and its tangent lines.

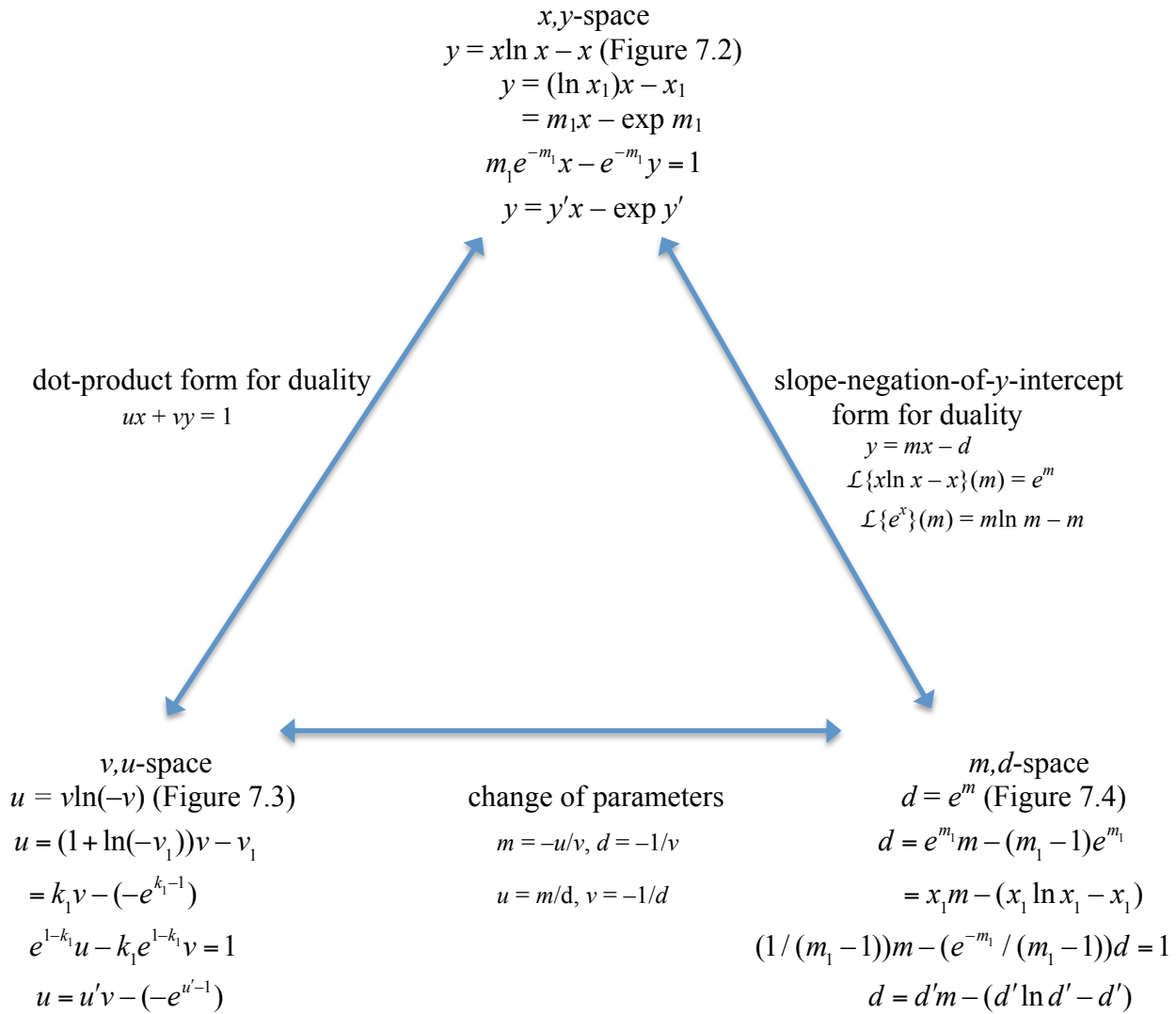

FIGURE 7.1: The curve $y = x\ln x - x$ and two of its dual curves

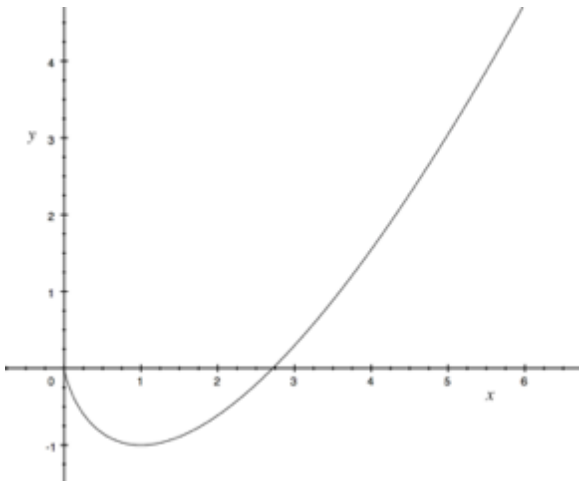

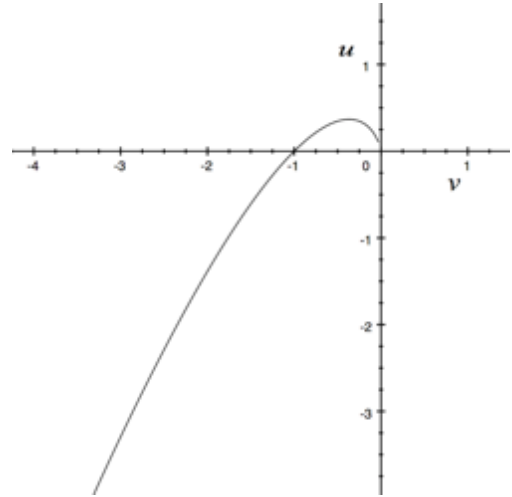

FIGURE 7.2: $y = x\ln x - x$ in $x,y$-space

FIGURE 7.3: $u = v\ln(-v)$ in $v,u$-space

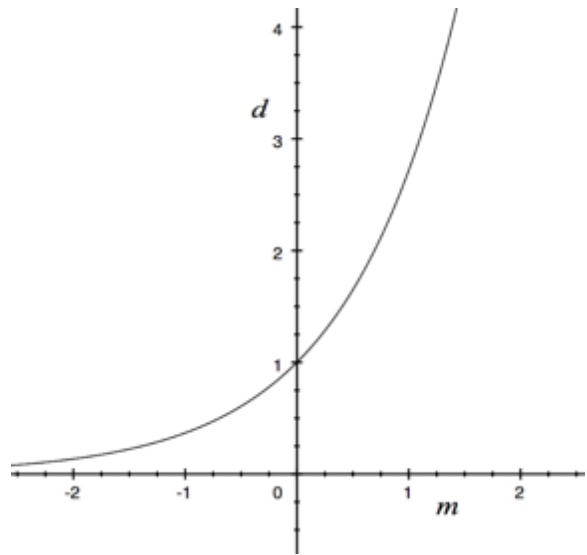

FIGURE 7.4: $d = e^m$ in $m,d$-space

## 8. Geometrically Finding the Coordinates of a Dual Curve

This section contains a geometric construction for points in *x,y*-space, whose coordinates are the same as the dual points for lines in *x,y*-space. The construction is a property of conics. A different conic is used for each duality relationship. This is explored for slope-negation-of-the-*y*-intercept duality and for dot-product duality. For the first duality relationship, in Section 8.1, point $(x,y) = (m,d)$ is constructed for the line $y = mx - d$. For the second duality relationship, in Section 8.2, the point $(x,y) = (u,v)$ is constructed for any line that can be written in the form $ux + vy = 1$.

Some of these ideas are ancient. Apollonius studied them [2, pp. 225–243; 8, p. 71; 20, pp. 102–108]. Also, Theorems 8.1–8.3 were investigated or proven in [6, pp. 120–122; 7, pp.199–230; 8, pp. 110–115, 331–332; 21, pp. 238–239; 35, pp. 70–73, 98–99, 121–122, 150, 194–195; 37, pp. 209–210; 39, pp. 251–253; 42]. The relationship between points and lines in the construction is the same as between dual points and tangent lines for curves and envelopes. One difference is that all the points, lines, and curves are in one space.

**Theorem 8.1** (*a pole determines its unique polar*): Select a point and draw chords to a fixed, non-degenerate conic through the point, which is called the *pole*. At the intersection of each chord and the conic, draw tangent lines to the conic. Then, the intersections of the pairs of tangent lines associated with each chord are on a line, which is called the *polar*. In particular, consider the non-degenerate conic

$$Ax^2 + 2Bxy + Cy^2 + 2Dx + 2Ey + F = 0 \qquad (8.1)$$

For each pole

$$(a,b), \qquad (8.2)$$

the polar with respect to the conic (8.1) is

$$\alpha x + \beta y + \delta = 0, \qquad (8.3)$$

where

$$\begin{pmatrix} \alpha \\ \beta \\ \delta \end{pmatrix} = \begin{pmatrix} A & B & D \\ B & C & E \\ D & E & F \end{pmatrix} \begin{pmatrix} a \\ b \\ 1 \end{pmatrix}. \qquad (8.4)$$

For a fixed conic, the relationship between poles and polars is one-to-one.

**Theorem 8.2** (*a polar determines its unique pole*): Select a line, which is the polar, and from the points on the line draw pairs of tangent lines to a fixed, non-degenerate conic. Draw the conic's chords through the points of tangency of each pair of tangent lines. Then, all those chords contain one point, which is the pole. For the non-degenerate conic (8.1), the pole of the line

$$\alpha x + \beta y + \delta = 0 \qquad (8.5)$$

with respect to the conic (8.1) is

$$\left( \frac{a}{c}, \frac{b}{c} \right), \qquad (8.6)$$

where

$$\begin{pmatrix} a \\ b \\ c \end{pmatrix} = \begin{pmatrix} A & B & D \\ B & C & E \\ D & E & F \end{pmatrix}^{-1} \begin{pmatrix} \alpha \\ \beta \\ \delta \end{pmatrix}. \qquad (8.7)$$

For a fixed conic, the relationship between poles and polars is one-to-one.

Theorems 8.1 and 8.2 describe an inverse relationship, which can be seen by substitution of (8.4) into (8.7) and of (8.7) into (8.4).

**Theorem 8.3** (*locations of poles and polars with respect to the conic*): If a pole (8.2) is inside the conic (8.1), then the polar (8.3) is entirely outside the conic. If the pole is on the conic, the polar is the tangent line at the pole. If the pole is outside the conic, the polar intersects the conic in two points.

Consider a family of tangent lines to a curve to be the polars. Each polar yields a pole. The set of the poles is a curve, which is called the *polar curve* of the original curve. The equation of the polar curve is the same as the corresponding dual curve, but in different coordinates. The difference between the dual curve and the polar curve is location. The dual curve is in the dual space, and the polar curve is in $x,y$-space.

## 8.1 The Pair: Pole $(m,d)$ and Polar $y = mx - d$

Consider the slope-negation-of-the-$y$-intercept standard form of the lines, which is the form for the dual space that is introduced in Section 3. The pole is $(m,d)$, and the polar is $y = mx - d$ or $mx - y - d = 0$. As noted in Table 7.1, vertical lines are excluded, because they do not possess a slope $m$. To find the coefficients of the appropriate conic (8.1) for the pole-polar correspondence of Theorems 8.1 and 8.2, substitute $(\alpha \ \beta \ \delta)^{-1} = (m \ -1 \ -d)^{-1}$ from the polar into the left-hand side of (8.4) and $(a \ b \ 1)^{-1} = (m \ d \ 1)^{-1}$ from the pole into the right-hand side of (8.4), giving

$$\begin{pmatrix} m \\ -1 \\ -d \end{pmatrix} = \begin{pmatrix} A & B & D \\ B & C & E \\ D & E & F \end{pmatrix} \begin{pmatrix} m \\ d \\ 1 \end{pmatrix},$$

which is

$$\begin{pmatrix} m \\ -1 \\ -d \end{pmatrix} = \begin{pmatrix} Am + Bd + D \\ Bm + Cd + E \\ Dm + Ed + F \end{pmatrix}.$$

This is an identity in $m$ and $d$, so that $A = 1$, $E = -1$, and $B = C = D = F = 0$. The conic (8.1) is
$$y = x^2/2.$$

A shortcut for finding the conic uses Theorem 8.3, which says that, when the polar is tangent to the conic, its pole is the point of tangency. Thus, if such a polar is $y = mx - d$, then its pole is $(m,d)$, which is on the polar. Substituting the coordinates of the pole into the equation of the polar gives $d = (m)(m) - d$, so that $d = m^2/2$. Because this is true for all points $(m,d)$ on the conic, the equation of the conic is $y = x^2/2$.

Say that we already know that the conic is $y = x^2/2$ and want the pole for the line $y = mx - d$. Then, (8.5) gives
$$(\alpha\ \beta\ \delta)^{-1} = (m\ -1\ -d)^{-1},$$
and (8.7) is
$$\begin{pmatrix} a \\ b \\ c \end{pmatrix} = \begin{pmatrix} 1 & 0 & 0 \\ 0 & 0 & -1 \\ 0 & -1 & 0 \end{pmatrix}^{-1} \begin{pmatrix} m \\ -1 \\ -d \end{pmatrix} = \begin{pmatrix} 1 & 0 & 0 \\ 0 & 0 & -1 \\ 0 & -1 & 0 \end{pmatrix} \begin{pmatrix} m \\ -1 \\ -d \end{pmatrix} = \begin{pmatrix} m \\ d \\ 1 \end{pmatrix},$$
so that the pole (8.6) is
$$\left(\frac{a}{c}, \frac{b}{c}\right) = (m,d).$$

The construction in Theorem 8.1 is illustrated in Figure 8.1. Line $L_1$ has the equation
$$y = mx - d. \tag{8.8}$$
The conic is the parabola
$$y = x^2/2. \tag{8.9}$$
The construction is given for the point $T(m,d)$ in $x,y$-space, which is the pole of the polar $L_1$ in (8.8) and has the same coordinates as the dual point $(m,d)$ in the dual space for slope-negation-of-the-$y$-intercept duality.

The line $L_2$ containing the two points $A(a,a^2/2)$ and $B(b,b^2/2)$ of the parabola (8.9) has the equation
$$y = ((a + b)/2)x - ab/2. \tag{8.10}$$
The tangent lines $L_3$ and $L_4$ at $A$ and $B$ are
$$y = ax - a^2/2 \text{ and } y = bx - b^2/2,$$
which meet at the point
$$R((a + b)/2, ab/2).$$
For the construction, $R$ is on the line (8.8), so that
$$ab/2 = m((a + b)/2) - d$$
and
$$b = (ma - 2d)/(a - m). \tag{8.11}$$
Substituting (8.11) into (8.10) gives
$$y = (1/2)((a^2 - 2d)/(a - m))x - a(ma - 2d)/(a - m), \tag{8.12}$$

which is a parameterization with parameter $a$ for fixed $m$ and $d$ of the lines $L_2$, which contain the chords in the construction of the pole. Placing the $x$-coordinate of the pole $T$, which is $m$, into (8.12) gives $y = d$. Thus, the point $T(m,d)$ is on all lines (8.12) that contain the chords and is the pole of $y = mx - d$ with respect to the parabola $y = x^2/2$. Two chords are sufficient for finding the intersection point $T$.

Alternatively, the pole for line (8.8) can be found as follows. Because $L_1$ in Figure 8.1 has $y$-intercept $-d$, the point $(0,d)$ can be found by reflecting $(0,-d)$ through the origin. The tangent of the angle that the ray through $(0,0)$ and $(m,d)$ makes with the positive $x$-axis is $d/m$. The intersection of $y = d$ and the line $y = (d/m)x$ is the point $(m,d)$.

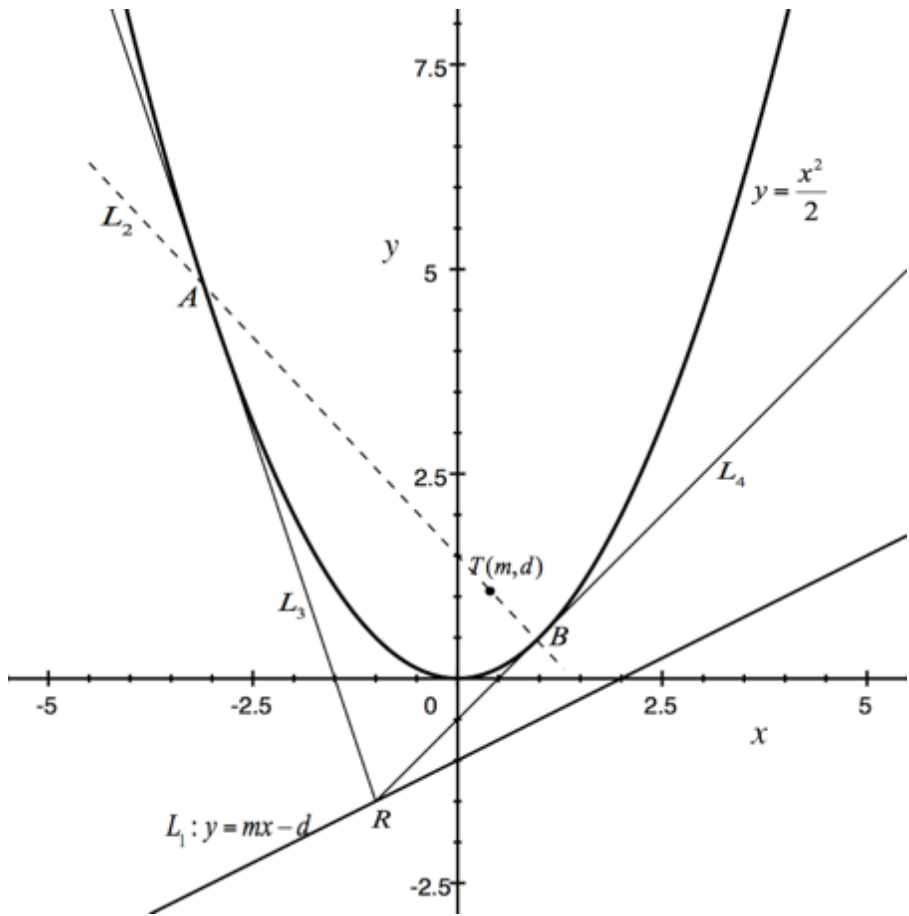

FIGURE 8.1: Polar $L_1$: $y = mx - d$ and pole $T(m,d)$ with respect to the conic $y = x^2/2$

## 8.2 The Pair: Pole $(u,v)$ and Polar $ux + vy = 1$

For the dot-product standard form of the lines, which is Entry 3 of Table 7.1, the pole is $(u,v)$ and the polar is $ux + vy = 1$. As noted in Table 7.1, lines containing the origin are excluded, because the $y$-intercept of $ux + vy = 1$ is $1/v \neq 0$. To find the coefficients of the conic (8.1) for the pole-polar correspondence of Theorems 8.1 and 8.2, substitute $(\alpha\ \beta\ \delta)^{-1} = (u\ v\ -1)^{-1}$ into the left-hand side of (8.4) and $(a\ b\ 1)^{-1} = (u\ v\ 1)^{-1}$ into the right-hand side of (8.4), giving

$$\begin{pmatrix} u \\ v \\ -1 \end{pmatrix} = \begin{pmatrix} A & B & D \\ B & C & E \\ D & E & F \end{pmatrix} \begin{pmatrix} u \\ v \\ 1 \end{pmatrix},$$

which is

$$\begin{pmatrix} u \\ v \\ -1 \end{pmatrix} = \begin{pmatrix} Au + Bv + D \\ Bu + Cv + E \\ Du + Ev + F \end{pmatrix}.$$

Because this is an identity in $u$ and $v$, $A = C = 1$, $F = -1$, and $B = D = E = 0$. The conic (8.1) is

$$x^2 + y^2 = 1.$$

As in Section 8.1, a shortcut for finding the conic uses Theorem 8.3, which says that, when the polar is tangent to the conic, its pole is the point of tangency. Thus, if such a polar is $ux + vy = 1$, then its pole is $(u,v)$, which is on the polar. Substituting the coordinates of the pole into the equation of the polar gives $(u)(u) + (v)(v) = 1$, so that $u^2 + v^2 = 1$. Because this is true for all points $(u,v)$ on the conic, the equation of the conic is $x^2 + y^2 = 1$.

Say that we already know that the conic is $x^2 + y^2 = 1$ and want the pole for the line $ux + vy = 1$. Then, (8.5) gives

$$(\alpha\ \beta\ \delta)^{-1} = (u\ v\ -1)^{-1},$$

and (8.7) is

$$\begin{pmatrix} a \\ b \\ c \end{pmatrix} = \begin{pmatrix} 1 & 0 & 0 \\ 0 & 1 & 0 \\ 0 & 0 & -1 \end{pmatrix}^{-1} \begin{pmatrix} u \\ v \\ -1 \end{pmatrix} = \begin{pmatrix} 1 & 0 & 0 \\ 0 & 1 & 0 \\ 0 & 0 & -1 \end{pmatrix} \begin{pmatrix} u \\ v \\ -1 \end{pmatrix} = \begin{pmatrix} u \\ v \\ 1 \end{pmatrix},$$

so that the pole (8.6) is

$$\left(\frac{a}{c}, \frac{b}{c}\right) = (u,v).$$

An alternative and simpler construction for $T(u,v)$ is in Theorem 8.4 [1, p. 35; 6, pp. 114–122; 16, p. 265; 27, p. 141; 35, pp. 72–73; 42]. Refer to Figure 8.2.

**Theorem 8.4** (*construction using inversion in the unit circle*): Consider the polar
$$L_1: ux + vy = 1.$$
Draw the line
$$L_2: y = (v/u)x,$$
which is perpendicular to the polar and contains the origin (0,0). The lines $L_1$ and $L_2$ meet at the point $Q$. Inverting $Q$ with respect to the unit circle yields the pole $T(u,v)$.

**Proof:** The coordinates of point $Q$ are $(u/(u^2 + v^2), v/(u^2 + v^2))$. The distance between (0,0) and $Q$ is $\dfrac{1}{\sqrt{u^2+v^2}}$. The distance between (0,0) and $T(u,v)$ is $\sqrt{u^2+v^2}$. The product of these distances is 1. ∎

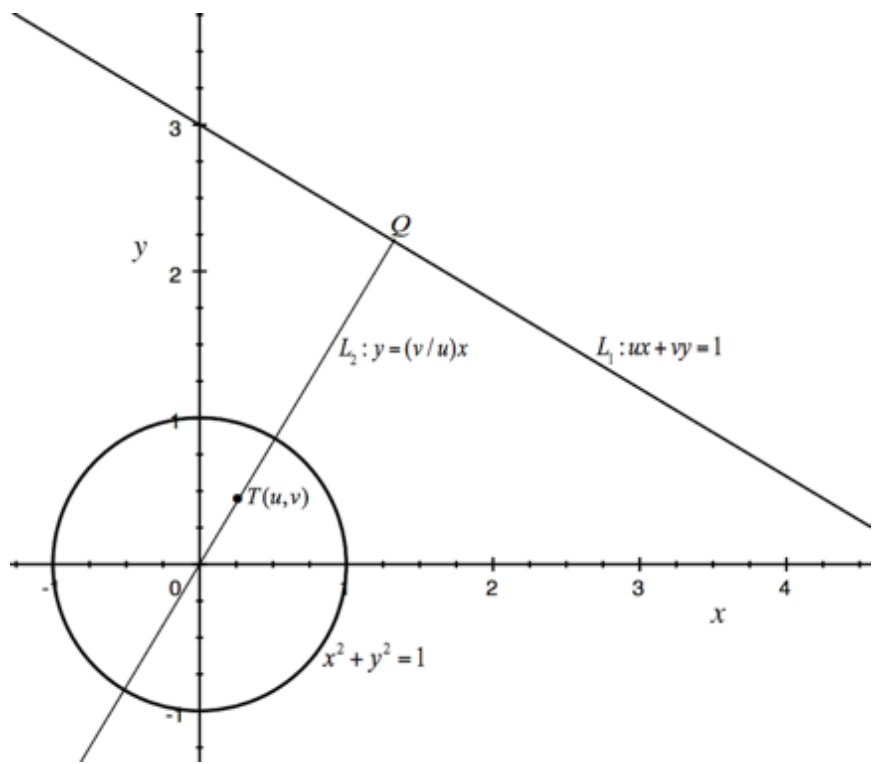

FIGURE 8.2: Alternative method to find the pole $T(u,v)$ to the polar
$L_1: ux + vy = 1$ with respect to the conic $x^2 + y^2 = 1$

## 9. Inequalities

The following four examples illustrate methods for finding and validating Legendre transformation pairs. The methods are founded upon the locations where inequalities are tight, that is, become equalities that are in the form of (3.1), which is $y(x) + d(m) = mx$, and hence, give pairs of dual functions. Theorem 9.1 justifies a method for finding pairs of dual functions by employing Young's inequality.

**Example 9.1** (*using the weighted arithmetic mean- geometric mean inequality*): Entry 7 of Table 5.1 is $\mathcal{L}\{x^p/p\}(m) = m^q/q$. We recognize the two functions $x^p/p$ and $m^q/q$ as terms in the inequality

$$\frac{x^p}{p} + \frac{m^q}{q} \geq xm \tag{9.1}$$

for $x > 0$, $m > 0$, $p > 1$, $q > 1$, and

$$1/p + 1/q = 1. \tag{9.2}$$

Inequality (9.1) can be derived from the weighted arithmetic mean – geometric mean inequality

$$r^\alpha s^\beta \leq \alpha r + \beta s \tag{9.3}$$

for positive $r$, $s$, $\alpha$, and $\beta$ and $\alpha + \beta = 1$ [19, pp. 37, 61; 30, pp. 27–30]. Setting $r = x^p$, $s = m^q$, $\alpha = 1/p$, and $\beta = 1/q$ in (9.3) gives (9.1) and gives (9.2) from $\alpha + \beta = 1/p + 1/q =1$. Because $0 < \alpha < 1$ and $0 < \beta < 1$, $p > 1$ and $q > 1$. There is equality in (9.3) if and only if $r = s$, so there is equality in (9.1) if and only if $x^p = m^q$, which is $m = x^{p/q}$. This means that $m$ is the slope of $x^p/p$, from

$$\frac{d}{dx}\frac{x^p}{p} = x^{p-1} = x^{p/q} = m.$$

**Example 9.2** (*finding the minimum of a function*): Another approach to (9.1) that appears in [19, p. 107], which yields $\mathcal{L}\{x^p/p\}(m) = m^q/q$ is to consider

$$f(x,m) = \frac{x^p}{p} + \frac{m^q}{q} - xm.$$

Because

$$\frac{\partial f(x,m)}{\partial x} = x^{p-1} - m \text{ and } \frac{\partial^2 f(x,m)}{\partial x^2} = (p-1)x^{p-2},$$

$\frac{\partial^2 f}{\partial x^2} > 0$ and the minimum value of zero for $f$ occurs for $\frac{\partial f}{\partial x} = 0$ [10, pp. 254–256; 12, pp. 158–163; 25, p. 194]. Thus, $m = x^{p-1}$ and

$$f(x,m) = \frac{x^p}{p} + \frac{m^q}{q} - xm = \frac{x^p}{p} + \frac{(x^{p-1})^q}{q} - xx^{p-1} = (\frac{1}{p}+\frac{1}{q})x^p - x^p = 0.$$

Examining $f(x,m)$ similarly as a function of $m$ further validates this conclusion.

**Example 9.3** (*finding the minimum of a function*): For Entry 13, which is $\mathcal{L}\{x\ln x\}(m) = e^{m-1}$,

$$f(x,m) = x\ln x + e^{m-1} - mx, \quad \frac{\partial f(x,m)}{\partial x} = \ln x + 1 - m, \text{ and } \frac{\partial^2 f(x,m)}{\partial x^2} = \frac{1}{x} > 0.$$

The minimum value of zero for $f$ occurs for $m = \ln x + 1$, that is,

$$f(x, \ln x + 1) = x\ln x + e^{\ln x + 1 - 1} - (\ln x + 1)x = 0$$

and

$$\frac{d}{dx} x\ln x = \ln x + 1 = m.$$

*Young's inequality* says that for the invertible, strictly increasing, and differentiable function $f$, which is defined on the interval $[0,c]$ with $f(0) = 0$, $0 < a < c$, and $0 < m < f(c)$,

$$\int_0^a f(t)\,dt + \int_0^m f^{-1}(t)\,dt \geq am \tag{9.4}$$

with equality if and only if

$$m = f(a) \tag{9.5}$$

[19, pp. 111–113; 30, pp. 48–50; 40, p. 14].

**Theorem 9.1** (*Young's inequality supplies dual functions*): Substituting for $f$ as described in Young's inequality into (9.4) with (9.5) yields an equation for a pair of dual functions under the Legendre transformation.

**Proof:** Define

$$y(x) = \int_0^x f(t)\,dt \tag{9.6}$$

for $0 \leq x \leq c$. From (9.6) for arbitrary $a \in [0,c]$,

$$y(a) = \int_0^a f(t)\,dt, \tag{9.7}$$

$$y'(x) = f(x), \tag{9.8}$$

and

$$y''(x) = f'(x) > 0.$$

From (9.8),

$$f^{-1}(t) = y'^{-1}(t).$$

Integrating gives

$$\int_0^m f^{-1}(t)\,dt = \int_0^m y'^{-1}(t)\,dt. \tag{9.9}$$

Substituting (9.7) and (9.9) into (9.4) and selecting $m$ according to (9.5), so that (9.4) is an equation, give

$$y(a) + \int_0^m y'^{-1}(t)\,dt = am. \qquad (9.10)$$

From (9.5) and (9.8), $m = f(a) = y'(a)$. From (5.5), the integral in (9.10) is

$$\mathcal{L}\{y(x)\}(m) = \int_0^m f^{-1}(t)\,dt = \int_0^m y'^{-1}(t)\,dt = d(m). \qquad (9.11)$$

Because $a$ is arbitrary, replace it with $x$. Combining (9.10) and (9.11) gives

$$y(x) + d(m) = mx$$

or

$$y(x) = mx - d(m). \blacksquare$$

The integral form of the Legendre transformation in (9.11) is also given in (5.5) in Corollary 5.1.

**Example 9.4** (*integral form of* $\mathcal{L}\{y(x)\}(m)$)**:** The strictly increasing function $f(t) = t^{p-1}$ with $p > 1$ and $t > 0$ gives Entry 7 of Table 5.1 with $x > 0$ and $m > 0$ by substitution into (9.6) and (9.11), so that

$$y(x) = \int_0^x f(t)\,dt = \int_o^x t^{p-1}\,dt = \frac{x^p}{p}$$

and

$$d(m) = \int_0^m f^{-1}(t)\,dt = \int_0^m t^{1/(p-1)}\,dt = \int_0^m t^{q-1}\,dt = \frac{m^q}{q}.$$

## 10. Supporting Lines and Convex and Concave Functions

An extension of the concept of a tangent line is a line of support. Using those lines allows us to expand our perspective to curves that have linear portions and corners as well as functions that are convex, but not differentiable.

### 10.1 Supporting Lines

A generalization of a tangent line is a supporting line to a curve. A *supporting line, support line,* or *line of support* is a line that contains a linear portion of a curve or else locally intersects a curve at just one point, while locally lying outside the curve. Tangent lines are supporting lines.

Curves with corners can be accommodated, because, even though there is no tangent line at the corner, there is a fan of supporting lines there [6, pp. 37–43; 15, pp. 20–21; 27, pp. 41–43, 140–146, 205–212; 38, p. 34]. The notion of a dual curve is extended to functions that have linear portions and corners by replacing tangent lines by supporting lines, as necessary, and employing the supporting lines and their coefficients analogously to tangent lines and their coefficients. Example 10.1 contains a piecewise linear function with three corners, and Example 10.2 is a piecewise linear curve with four corners. The examples illustrate how linear portions (corner points) of a curve in $x,y$-space correspond to corner points (linear portions) of a curve in $m,d$-space.

**Example 10.1** (*supporting lines give the dual curve to a piecewise linear function*)**:** Consider the piecewise linear function that is defined by

$$y = y(x) = \begin{cases} -2x - 2 & x \leq -3 \\ -x + 1 & -3 < x \leq 0 \\ x + 1 & 0 < x \leq 1 \\ 3x - 1 & 1 < x \end{cases} \quad (10.1)$$

and appears in Figure 10.1.

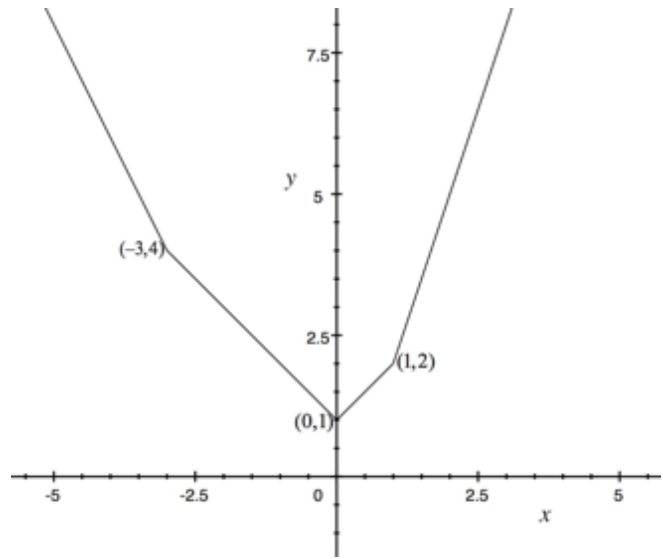

FIGURE 10.1: The piecewise-linear function (10.1)

The dual points of the four sides are $(m,d) = (-2,2)$ for $y = mx - d = -2x - 2$, $(-1,-1)$ for $y = -x + 1$, $(1,-1)$ for $y = x + 1$, and $(3,1)$ for $y = 3x - 1$. The supporting lines at $(x,y) = (-3,4)$ are $y = m(x-(-3)) + 4 = mx + 3m + 4$ for $-2 < m < -1$, whose dual points are on the line segment with $d = -3m - 4$. The supporting lines at $(0,1)$ are $y = mx + 1$ for $-1 < m < 1$, whose dual points are on the line segment with $d = -1$. The supporting lines at $(1,2)$ are $y = mx - m + 2$, whose dual points are on the line segment with $d = m - 2$. The dual curve, which has domain $-2 \leq m \leq 3$ and appears in Figure 10.2, is

$$d = d(m) = \begin{cases} -3m - 4 & -2 \leq m \leq -1 \\ -1 & -1 < m \leq 1 \\ m - 2 & 1 < m \leq 3 \end{cases} \quad (10.2)$$

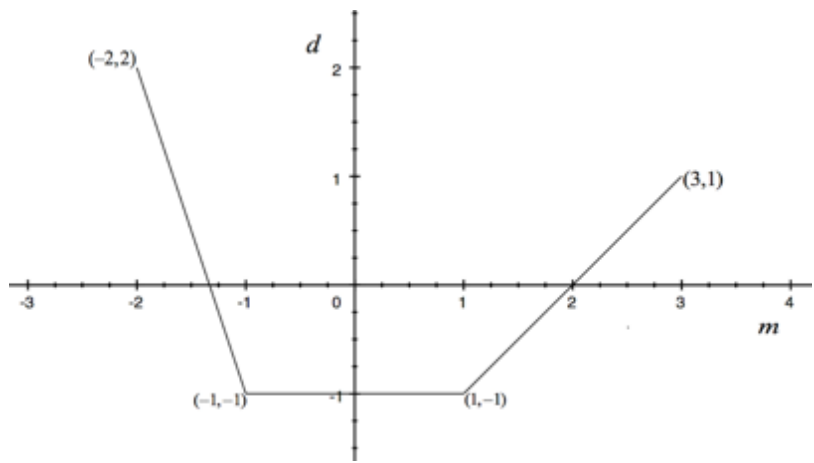

FIGURE 10.2: The dual function (10.2)

**Example 10.2** (*supporting lines give the dual curve to a piecewise linear curve*)**:** An example of a piecewise linear curve appears in Figure 10.3. The line $y = x + 1$, which contains the upper-left side for $-1 < x < 0$, has the dual point $(m,d) = (1,-1)$, which is the lower-right vertex in the dual curve in Figure 10.4. The supporting lines through the vertex $(x,y) = (0,1)$ have equations $y = mx + 1$ for $-1 < m < 1$, whose dual points are $(m,d) = (m, -1)$ for $-1 < m < 1$, which is the bottom side containing $(-1,-1)$ and $(1,-1)$ in Figure 10.4.

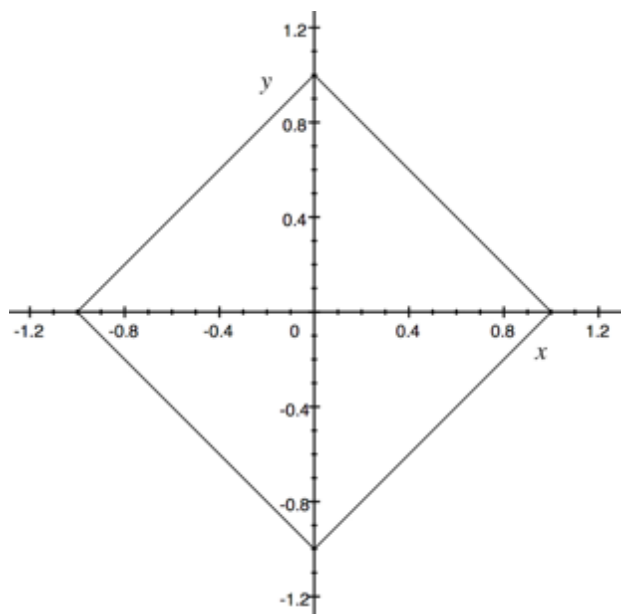

FIGURE 10.3: A piecewise-linear curve

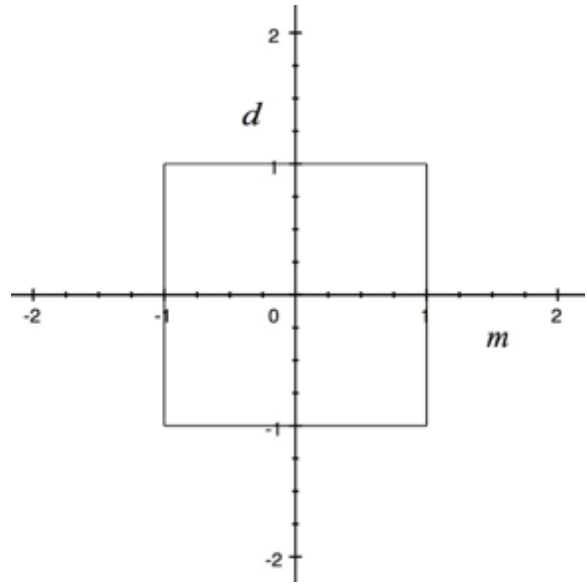

FIGURE 10.4: The dual curve of the curve in Figure 10.3

## 10.2 Convex and Concave Functions

Examples 10.1 and 10.2 demonstrate the creation of a dual curve by using lines of support to a curve that is not everywhere differentiable. Example 10.3 and the revisited Example 10.1 illustrate a different approach for finding a dual function, which can be applied to convex or concave non-differentiable functions. For a convex function and each value of slope $m$, a sequence of lines that are below the function are drawn (or imagined). The limiting line which has the largest negation of the $y$-intercept produces the value of $-d(m)$. This process yields the envelopes of convex functions. Convex functions are of necessity continuous, and there exists a line of support at each of its points [5, p. 151; 6, pp. 37–43; 15, p. 27; 27, pp. 41–42; 34, pp. 102–103; 40, pp. 11, 17]. Because concave functions can be expressed as the negation of a convex function, they too have supporting lines and hence envelopes, as well. The methodology, which is illustrated in the two examples below, can be used over each interval where a curve is concave or convex.

**Example 10.3** (*the negation of the y-intercept at the limit of a sequence of lines below a convex function*)**:** For illustrative purposes and simplicity, consider the function
$$y = y(x) = x^3/3 \text{ for } 0 \leq x \leq 2,$$
which is differentiable and displayed in Figure 10.5. From Entry 7 of Table 7.1 with $p = 3$ and $q = 3/2$, the dual curve is
$$d = d(m) = m^{3/2}/(3/2) \text{ for } 0 \leq m \leq 4, \qquad (10.3)$$
where the bounds on the $m$-interval are the slopes of $y = y(x)$ at its endpoints (0,0) and (2,8/3).

Select a value for *m* between 0 and 4 and draw a succession of lines with slope *m* that are below or just touching the curve. In Figure 10.5, $m = 1$ and four lines are drawn. As the *y*-intercept of those lines increases, the lines approach a supporting line, which is a tangent line in this case. In the limit, the point of contact (1,1/3) is reached, the *y*-intercept for that line is –2/3, and $d(1) = 2/3$, which is substantiated by (10.3).

When this approach is used, the domain of $d = d(m)$ is extended to the real line by using the lines of support at the end points of the domain of $y = y(x)$, which is a finite, closed *x*-interval here. The dual curve is

$$d = d(m) = \begin{cases} 0 & m < 0 \\ \dfrac{2m^{3/2}}{3} & 0 \le m \le 4 \\ 2m - \dfrac{8}{3} & 4 < m \end{cases} \qquad (10.4)$$

All lines through $(x,y) = (0,0)$ with $m < 0$ yield $d = d(m) = 0$, and lines through $(x,y) = (2,8/3)$ with $m > 4$ have *y*-intercepts $8/3 - 2m$ and thus $d = d(m) = 2m - 8/3$. Function (10.4) is continuous and convex. It is convex because $y = y(x)$ is convex, as implied by Theorem 3.1.

The important feature of $y = y(x)$ for this construction is that it is convex over the whole domain of interest. That assures that the sequence of lines with the given slope has a limiting line from below.

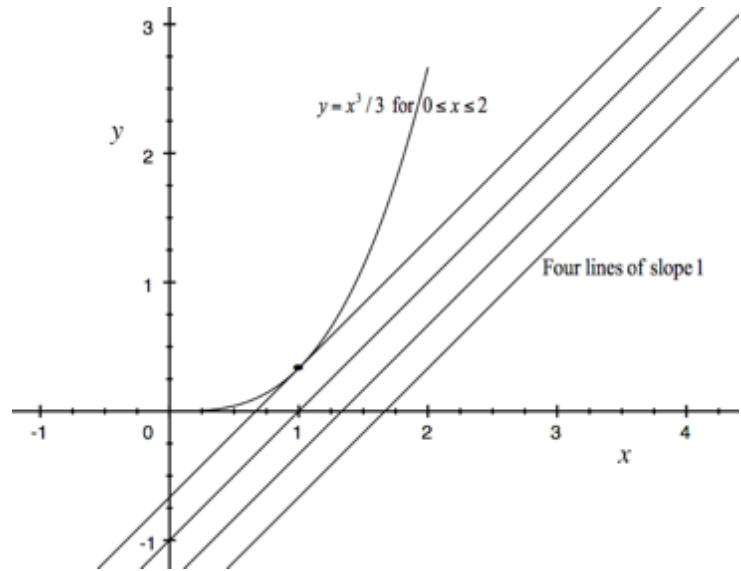

FIGURE 10.5: The curve $y = x^3/3$ for $0 \le x \le 2$ and four lines with slope 1 that do not cross the curve

Continuing this example and writing $y = mx - g$, the negation of the *y*-intercept for a given *m* with $0 < m < 4$ is

$$g(x) = mx - x^3/3.$$

To find its maximum value, compute

$$g'(x) = \frac{d}{dx}(mx - y(x)) = m - y'(x) = m - x^2 = 0,$$

which is the familiar $y'(x) = m$ and $x = m^{1/2}$. The second derivative $g''(x) = -2x < 0$ assures that $g(x)$ is a concave function, whose maximum is reached at $x = m^{1/2}$, yielding $d(m) = m(m^{1/2}) - (m^{1/2})^3/3 = 2m^{3/2}/3$, which is confirmed by (10.3).

**Example 10.1 continued** (*the negation of the y-intercept at the limit of a sequence of lines below a convex, piecewise-linear function*): A sequence of lines with slope $m = 5/2$ and values less than or equal to convex function (10.1) in Figure 10.1 has a member that intersects the function at (1,2). The y-intercept of that line from $y = mx + b$ is $b = 2 - (5/2)(1) = -1/2$. Hence, $d(5/2) = 1/2$, which is substantiated by (10.2).

The negation of the y-intercept of lines of slope $m$ that lie below the convex function (10.1), which appears in Figure 10.1, is

$$g = mx - y(x) = \begin{cases} (m+2)x + 2 & x \le -3 \\ (m+1)x - 1 & -3 < x \le 0 \\ (m-1)x - 1 & 0 < x \le 1 \\ (m-3)x + 1 & 1 < x \end{cases} \quad (10.5)$$

For each value of $m$, the maximum value of $g(x) = mx - y(x)$ is the value of the function $d(m)$ that is dual to $y = y(x)$. For $m = 5/2$, (10.5) is

$$g = \frac{5}{2}x - y(x) = \begin{cases} \dfrac{9}{2}x + 2 & x \le -3 \\ \dfrac{7}{2}x - 1 & -3 < x \le 0 \\ \dfrac{3}{2}x - 1 & 0 < x \le 1 \\ -\dfrac{1}{2}x + 1 & 1 < x \end{cases} \quad (10.6)$$

The maximum of (10.6) occurs at $x = 1$, so that $d(5/2) = 1/2$.

For the convex function $y = y(x)$, the notation is

$$d(m) = \sup_x \{mx - y(x)\}$$

where "sup" indicates "supremum," which is also called the "least upper bound" or "lub" [6, p. 7; 15, p. 49]. This gives an alternative point-of-view of this duality [3, pp. 61–62; 5, p. 219; 41]. For each value of $m$, $d(m)$ is the least upper bound of all the vertical signed distances from $y = mx$ to $y = y(x)$. For differentiable functions $y(x)$, $d/dx(mx - y(x)) = 0$ yields the x-value for the maximum and $x = y'^{-1}(m)$. Thus, the same function $d(m)$ is obtained. Figure 10.6 illustrates that the two definitions of $d(m)$, that is, in Section 3 as the negation of the y-intercept and as the supremum of the distances between $y = mx$ to $y = y(x)$, give the same valuations, but are measured differently. They are the lengths of opposite sides of a parallelogram.

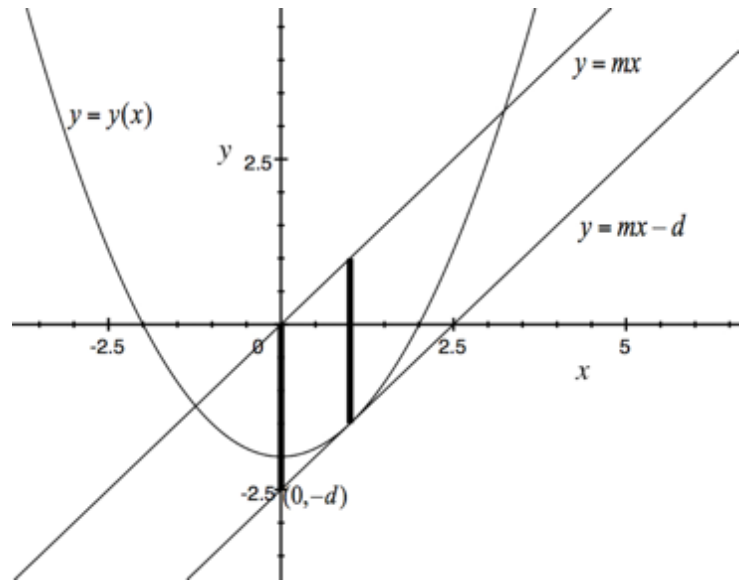

FIGURE 10.6: The two line segments in bold give two equal measurements for $d(m)$
and thus yield the same value for $d(m)$

This formulation and demonstration with Figure 10.6 make the proofs of some results almost trivial. An example is the next theorem [5, p. 232].

**Theorem 10.1** (*reversal of inequality between convex functions in the two spaces*): For convex functions $f$ and $g$, with $f(x) \leq g(x)$ for all values in their common domain if and only if $\mathcal{L}\{f(x)\}(m) \geq \mathcal{L}\{g(x)\}(m)$ for all $m$ in the common domain.

For the strictly concave function $y = y(x)$, write
$$\mathcal{L}\{y(x)\}(m) = d(m) = \inf_x \{mx - y(x)\}$$
where "inf" indicates "infimum," which is also called the "greatest lower bound" or "glb," and the analysis proceeds similarly.

## 11. Closing Comments

The nomenclature for many of these concepts is not fixed. For example, Bix [8, p. 110] calls the tangent lines the *envelope* of the curve. Occasionally, the function $b = b(m)$ is called the *Legendre transformation* of $y = y(x)$ [31]. To add further possible confusion, the name *Legendre transformation* is also used for the integral Legendre transform in which the kernels are Legendre polynomials [13, pp. 513–528]. Often, the *polar curve* is called the *dual curve* [1, pp. 68–69; 8, pp. 308–310] and the *polar reciprocal* [7, p. 217]. Eggleston [15, p. 20] and Benson [6, p. 114] use *polar duality* for duality. Van Tiel [40, p. 84] takes the words conjugate, dual, and polar function as synonymous. In *m,d*-duality, alternative symbols are commonly used for *m* (*e.g.*, $x^*$ and $p$) and for *d* (*e.g.*, $y^*$, $^*y$, and $f^*$). The variables $x$ and $m$ are sometimes called *conjugate variables*. There should be no misunderstanding, because all those cited works contain clearly stated definitions.

The main theme of Sections 1–10 is the connections among the ideas of family of tangent and support lines, envelope, dual space, dual function, differential equation, polar curve, and Young's inequality made tight. Study could be focused on the slope-negative-*y*-intercept form for lines and the Legendre transformation, because Artstein-Avidan and Milman [4] showed the centrality of the Legendre transformation between curves under very general conditions, and there are transformations between the various dual spaces, as shown in Section 7.

Another theme or thread is a variety of methods for finding the dual of a function or a pair of dual functions in *x,y*- and *m,d*-spaces. There is the straightforward expression

$$\mathcal{L}\{y(x)\}(m) = d(m) = mx - y = m{y'}^{-1}(m) - y({y'}^{-1}(m))$$

in Theorem 3.2, the integral in Corollary 5.1, dual functions that are obtained from known functions using Entries 1–6 of Table 5.1, solutions to Clairaut's differential equation from Theorem 6.1, transforming a pair of dual curves from another dual space as discussed in Section 7, the geometric construction in Theorem 8.2 and Section 8.1 using the parabola $y = x^2/2$, substitution of a function into the integrals in Young's inequality as described in Theorem 9.1 and Example 9.1, and the limiting procedure in Section 10.2.

Almost all the ideas here can be applied in *n*-dimensions to surfaces and tangent and supporting hyperplanes [6; 15; 24; 27; 34; 38; 40; 43]. The Lie theory of one-parameter continuous groups of symmetries for differential equations has yielded many insights [11; 22; 23; 32], including into Clairaut's differential equation [9], and might be successfully applied to the subject of duality. It is natural to consider these theorems in the projective plane [1; 8; 16; 35, pp. 315–338].

**Acknowledgement**

The authors wish to express their gratitude and appreciation to Patricia M. Burgess for her many insightful comments and suggestions on drafts of this paper.**References**

1.  A. V. Akopyan and A. A. Zaslavsky, *Geometry of Conics*, American Mathematical Society (2007).
2.  Apollonius of Perga, C. R. Talaferro, and M. N. Fried, *Conics: Books* I–IV, Green Lion (2013).
3.  V. I. Arnold, *Mathematical Methods of Classical Mechanics*, 2nd edn, Springer (1989).
4.  S. Artstein-Avidan and V. Milman, The Concept of Duality in Convex Analysis, and the Characterization of the Legendre Transform, *Annals of Mathematics, Second Series* 169 (2, March) (2009) 661–674.
5.  H. H. Bauschke and P. L. Combettes, *Convex Analysis and Monotone Operator Theory in Hilbert Space*, 2nd edn, Springer (2017).39

29. E. Maor, Line Equations of Curves: Duality in Analytic Geometry, *International Journal of Mathematics Education in Science and Technology* 9 (3) (1978) 311–322.

30. D. S. Mitrinović, *Analytic Inequalities*, Springer (1970).

31. C. E. Mungan, Legendre transforms for dummies (2014) 11 pages, available at www.aapt.org/docdirectory/meetingpresentations/SM14/Mungan-Poster.pdf (accessed December 5, 2020).

32. P. J. Olver, *Applications of Lie Groups to Differential Equations*, Springer (2000).

33. E. D. Rainville, *Elementary Differential Equations*, 3rd edn., Macmillan (1964).

34. R. T. Rockafellar, *Convex Analysis*, Princeton University Press (1970).

35. C. Smith, *Elementary Treatise on Conic Sections*, 2nd edn., Macmillan (1902).

36. S. Sternberg, Legendre Transformations of Curves, *Proceedings of the American Mathematical Society* 5 (6, December) (1954) 942–945.

37. J. H. Tanner and J. Allen, *Analytic Geometry*, American Book Company (1898).

38. A. C. Thompson, *Minkowski Geometry*, Cambridge (1996).

39. I. Todhunter, *A Treatise on Plane Co-Ordinate Geometry as Applied to the Straight Line and the Conic Sections*, 7th edn., Macmillan (1881).

40. J. van Tiel, *Convex Analysis*: *An Introductory Text*, Wiley (1984).

41. Wikipedia contributors, Legendre Transformation, *Wikipedia, The Free Encyclopedia,* https://en.wikipedia.org/w/index.php?title=Legendre_transformation&oldid=984539190 (accessed December 5, 2020).

42. Wikipedia contributors, Pole and Polar, *Wikipedia, The Free Encyclopedia,* https://en.wikipedia.org/w/index.php?title=Pole_and_polar&oldid=964110576 (accessed December 5, 2020).

43. R. K. P. Zia, E. F. Redish, and S. R. McKay, Making Sense of the Legendre Transform. *American Journal of Physics* 77 (July) (2009) 614–622.

44. D. G. Zill, *Differential Equations with Boundary-Value Problems*, Prindle, Weber, and Schmidt (1986).